\documentclass[preprint,10pt]{elsarticle}

\usepackage{url}
\usepackage{graphicx}
\usepackage{caption}
\usepackage{subfig}
\usepackage{epsfig}
\usepackage{psfrag}
\usepackage{hyperref}
\usepackage{amsmath,amssymb,amsfonts,amsthm,stmaryrd}
\usepackage{dsfont}
\usepackage{eucal}
\usepackage{mathrsfs}
\usepackage{listings}    
\usepackage{color}
\usepackage{float}
\usepackage{natbib}
\usepackage{booktabs}
\usepackage{tabularx}
\usepackage{setspace}
\usepackage{algorithmicx}
\usepackage{algorithm}
\usepackage{algpseudocode}

\makeatletter
\renewcommand\fs@ruled{\def\@fs@cfont{\bfseries}\let\@fs@capt\floatc@ruled
\def\@fs@pre{\hrule height 1.2pt depth0pt \kern2pt}%
\def\@fs@post{\kern2pt\hrule height 1.2pt depth0pt \kern2pt \relax}%
\def\@fs@mid{\kern2pt\hrule\kern2pt}%
\let\@fs@iftopcapt\iftrue}
\makeatother

\newcommand{\norm}[1]{\left|\left|#1\right|\right|}

\newcommand{\cref}[1]{Chapter~\ref{chapt:#1}}

\newcommand{\D}{\displaystyle}

\newcommand{\di}{\mathrm{d}}

\usepackage{xspace}
\newcommand{\eg}{e.g.,\xspace}

\newcommand{\ie}{i.e.,\xspace}

\newcommand{\cmatrixb}{\left\{ \begin{matrix}}
\newcommand{\cmatrixe}{\end{matrix} \right\}}

\newcommand{\vm}[1]{\mathbf{#1}}

\newcommand{\bsym}[1]{\boldsymbol{#1}}

\newcommand{\trans}{^\mathrm{T}}

\newcommand{\bc}{\begin{center}}
\newcommand{\ec}{\end{center}}
\newcommand{\bitem}{\begin{itemize}}
\newcommand{\eitem}{\end{itemize}}

\newcommand{\ljump}{\lbrack \! \lbrack } 
\newcommand{\rjump}{\rbrack \! \rbrack } 
\newcommand{\jump}[1]{\ljump {#1} \rjump} 

\newcommand{\beq}{\begin{equation}}
\newcommand{\eeq}{\end{equation}}
\newcommand{\beqa}{\begin{eqnarray}}
\newcommand{\eeqa}{\end{eqnarray}}

\newcommand{\invisible}[1]{}

\newcommand{\bv}{\begin{verbatim}}
\newcommand{\V}{\verb}                  % EX: \V=-d{#@~}= Expr must

\newcommand{\testpix}[1]{\fbox{\begin{minipage}[c]{\textwidth}
                      #1 \end{minipage} }}

\newcommand{\putpstex}[1]{\includegraphics{#1.pstex_t}}

\newcommand{\grbf}[1]{\mbox{\boldmath{$#1$}}}

\floatstyle{ruled}
\newfloat{Program}{thp}{lop}
\floatname{Program}{Program}

\floatstyle{ruled}
\newfloat{Fbox}{thp}{lop}
\floatname{Fbox}{Box}

\newcolumntype{C}{>{\centering\arraybackslash}X}

\definecolor{darkgray}{rgb}{0.95,0.95,0.95}
\definecolor{mygreen}{rgb}{0,0.6,0}

\lstdefinestyle{Matlab}
{
 basicstyle=\footnotesize, numbers=none, numberstyle=\tiny,%
 showstringspaces=false, language=Matlab, escapechar=|,frame=tb,%
 commentstyle=\color{mygreen}
}

\lstdefinestyle{Matlab-num}
{
 basicstyle=\footnotesize, numbers=left, numberstyle=\tiny,%
 showstringspaces=false, language=Matlab, escapechar=|,frame=tb,%
commentstyle=\color{mygreen}
}

\newcommand{\tty}[1]{\textnormal{\texttt{#1}}}
\newcommand{\sym}[1]{\textnormal{\textit{#1}}}

\lstnewenvironment{snippet}[1][]
{
 \lstset{style=Matlab, xleftmargin=5mm, gobble=4, #1}
}
{}

\lstnewenvironment{snippet1}[1][]
{
 \lstset{style=Matlab-num, xleftmargin=5mm, gobble=4, #1}
}
{}

\def\bibsection{\section*{References}}

\begin{document}


\definecolor{MyDarkBlue}{rgb}{1, 0.9, 1}
\lstset{language=Matlab,
       basicstyle=\footnotesize,
       commentstyle=\itshape,
       stringstyle=\ttfamily,
       showstringspaces=false,
       tabsize=2}
\lstdefinestyle{commentstyle}{color=\color{green}}

\theoremstyle{remark}
\newtheorem{thm}{Theorem}[section]
\newtheorem{rmk}[thm]{Remark}


\definecolor{red}{gray}{0}
\definecolor{blue}{gray}{0}


\begin{frontmatter}

\title{Nitsche's method for two and three dimensional NURBS patch coupling}       

\author[cardiff]{Vinh Phu Nguyen  \fnref{fn1}}
\author[cardiff]{Pierre Kerfriden  \fnref{fn2}}
\author[torino]{Marco Brino  \fnref{fn4}}
\author[cardiff]{St\'{e}phane P.A. Bordas \corref{cor1}\fnref{fn3}}
\author[torino]{Elvio Bonisoli  \fnref{fn5}}

\cortext[cor1]{Corresponding author}

\address[cardiff]{School of Engineering, Institute of Mechanics and Advanced
Materials, Cardiff University, Queen's Buildings, The Parade, Cardiff \\
CF24 3AA}

\address[torino]{Politecnico di Torino - DIGEP, corso Duca degli Abruzzi 24, 10129 Torino}

\fntext[fn1]{\url nguyenpv@cardiff.ac.uk, ORCID: 0000-0003-1212-8311}
\fntext[fn2]{\url pierre@cardiff.ac.uk}
\fntext[fn3]{\url stephane.bordas@alum.northwestern.edu, ORCID: 0000-0001-7622-2193}
\fntext[fn4]{\url marco.brino@polito.it}
\fntext[fn5]{\url elvio.bonisoli@polito.it}

\maketitle
\tableofcontents

\begin{abstract}
A Nitche's method is presented to couple non-conforming two and three dimensional NURBS (Non Uniform Rational B-splines)
patches in the context of isogeometric analysis (IGA). We present results for elastic stress analyses under the static condition
of two and three dimensional NURBS geometries. The contribution fills the gap in the literature and enlarges the applicability of
NURBS-based isogeometric analysis.
\end{abstract}

\begin{keyword} 
   Nitsche's method \sep isogeometric analysis (IGA) \sep multi-patch NURBS IGA \sep finite element method
\end{keyword}

\end{frontmatter}


\section{Introduction}

The predominant technology that is used by CAD to represent complex
geometries is the Non-Uniform Rational B-spline (NURBS). This allows certain
geometries to be represented exactly that are only approximated by
polynomial functions, including conic and circular sections. There is a
vast array of literature focused on NURBS
(e.g. \cite{piegl_book}, \cite{Rogers2001}) and as a result of several decades
of research, many efficient computer algorithms exist for their fast evaluation and
refinement. The key concept outlined by Hughes et al.
\cite{hughes_isogeometric_2005} was to employ NURBS not only as a geometry discretisation technology, but also as a discretisation tool for analysis, attributing such methods to the field of `Isogeometric Analysis'
(IGA). Since this seminal paper, a monograph
dedicated entirely to IGA has been
published \cite{cottrel_book_2009} and applications can now be found
in several fields including structural mechanics, solid
mechanics, fluid mechanics and contact mechanics. 
It should be emphasized that the idea of using CAD technologies in finite elements dates back
at least to \cite{NME:NME292,Kagan2000539} where B-splines were used as shape functions in FEM. In addition, similar methods which adopt
subdivision surfaces have been used to model shells \cite{Cirak_2000}. 

Structural mechanics is a field where IGA has demonstrated 
compelling benefits over conventional approaches
\cite{benson_isogeometric_2010,kiendl_isogeometric_2009,benson_large_2011,
beirao_da_veiga_isogeometric_2012,uhm_tspline_2009,Echter2013170,Benson2013133}. 
The smoothness of the NURBS basis functions allows for a straightforward
construction of plate/shell elements. Particularly for thin shells, rotation-free
formulations can be easily constructed \cite{kiendl_isogeometric_2009,kiendl_bending_2010}.
Furthermore, isogeometric plate/shell elements exhibit much less
pronounced shear-locking compared to standard FE plate/shell elements.

In  contact formulations using conventional geometry discretisations, the presence of 
faceted surfaces can lead to jumps and
oscillations in traction responses unless very fine meshes are used. The benefits of
using NURBS over such an approach are evident, since smooth contact
surface are obtained, leading to more physically accurate contact stresses. Recent work in this area includes 
\cite{temizer_contact_2011,jia_isogeometric_2011,temizer_three-dimensional_2012,
de_lorenzis_large_2011,Matzen201327}.

IGA has also shown advantages over traditional approaches in the context of optimisation problems 
\cite{wall_isogeometric_2008,manh_isogeometric_2011,qian_isogeometric_2011,xiaoping_full_2010} where the tight coupling with CAD
models offers an extremely attractive approach for industrial
applications.  Another attractive class of methods include those that require only a boundary discretisation, creating a truly direct coupling with CAD. Isogeometric boundary element methods for elastostatic analysis were presented in 
\cite{simpson_two-dimensional_2012,Scott2013197}, demonstrating that mesh generation can be completely circumvented by using CAD discretisations for analysis.

The smoothness of NURBS basis functions is attractive for analysis of 
fluids \cite{gomez_isogeometric_2010,nielsen_discretizations_2011,Bazilevs:2010:LES:1749635.1750210} and for
fluid-structure interaction problems 
\cite{bazilevs_isogeometric_2008,bazilevs_patient-specific_2009}. 
In addition, due to the ease of constructing high order continuous basis functions, IGA has been 
used with great success in solving PDEs that incorporate fourth order (or
higher)
derivatives of the field variable such as the Hill-Cahnard equation 
\cite{gomez_isogeometric_2008}, explicit gradient damage models \cite{verhoosel_isogeometric_2011-1} and gradient
elasticity \cite{fischer_isogeometric_2010}. 
The high order NURBS basis has also found potential applications in the Kohn-Sham equation for electronic 
structure modeling of semiconducting materials \cite{Masud2012112}.

NURBS provide advantageous properties for structural vibration problems 
\cite{cottrell_isogeometric_2006,Hughes20084104,NME:NME4282,Wang2013}
where $k$-refinement (unique to IGA)
has been shown to provide more robust and accurate frequency spectra than 
typical higher-order FE $p$-methods. Particularly, the optical branches of frequency spectra,
which have been identified as contributors to Gibbs phenomena in wave propagation problems 
(and the cause of rapid degradation of higher modes in the $p$-version of FEM), 
are eliminated. However when lumped mass matrices were used, the accuracy is limited to second order
for any basis order. High order isogeometric lumped mass matrices are not yet available.
The mathematical properties of IGA were studied in detail by Evans et al.\cite{evans_n-widths_2009}.

IGA has been applied to cohesive fracture \cite{verhoosel_isogeometric_2011}, outlining a framework for
modeling debonding along material interfaces using NURBS and propagating cohesive
cracks using T-splines. The method relies upon the ability to specify the continuity of NURBS and T-splines
through a process known as knot insertion.
As a variation of the eXtended Finite Element Method (XFEM)
\cite{mos_finite_1999}, IGA was applied to Linear Elastic Fracture Mechanics (LEFM) using the partition
of unity method (PUM) to capture two dimensional strong discontinuities and
crack tip singularities efficiently \cite{de_luycker_xfem_2011,ghorashi_extended_2012}. The method is usually
referred to as XIGA (eXtended IGA).
In \cite{Tambat20121} an explicit isogeometric enrichment technique was proposed for modeling
material interfaces and cracks exactly. Note that this method
is contrary to PUM-based enrichment methods which define cracks implicitly.
A phase field model for dynamic fracture was presented in \cite{Borden201277} using adaptive T-spline refinement to provide an effective method for simulating fracture in three dimensions.
In \cite{Nguyen2013} high order B-splines were adopted to efficiently 
model delamination of composite specimens and in \cite{nguyen_cohesive_2013}, an isogeometric framework
for two and three dimensional delamination analysis of composite laminates was presented where
the authors showed that using IGA can significantly reduce the usually time consuming pre-processing
step in generating FE meshes (solid elements and cohesive interface elements) for delamination computations.
A continuum description of fracture using 
explicit gradient damage models was also studied using NURBS \cite{verhoosel_isogeometric_2011-1}.

In computer aided geometric design, objects of complex topologies are usually represented as 
multiple-patch NURBS. We refer to Fig.~\ref{fig:concepts} for such a multi-patch NURBS solid.
Since it is virtually impossible to have a conforming parametrisation at the patch interface,
an important research topic within the IGA context is the implementation of multi-patch methods 
with high inter-patch continuity properties.
In this paper, a Nitsche's method is presented to couple non-conforming two and three dimensional NURBS patches
in a weak sense. An exact multipoint constraint method was reported in \cite{cottrel_book_2009} to glue
multiple NURBS patches with the restriction that, in the coarsest mesh, they have the same parametrisation.
Another solution to multi-patch IGA which has gathered momentum from both the computational geometry and analysis communities is the use of T-splines \cite{Sederberg:2003:TT:882262.882295}. T-splines correct the
deficiencies of NURBS by creating a single patch, watertight geometry which can be locally refined and coarsened. 
Utilisation of T-splines in an IGA framework has been illustrated in
\cite{bazilevs_isogeometric_2010,doerfel_adaptive_2010,scott_isogeometric_2011}. However T-splines are not yet a standard in CAD
and therefore our contribution will certainly enlarge the application areas of NURBS based IGA. Moreover, the formulation presented in
this contribution lays the foundation for the solid-structure coupling method to be presented in a forthcoming paper \cite{nguyen-nitsche2}.

\begin{figure}
  \centering
  \includegraphics[width=0.8\textwidth]{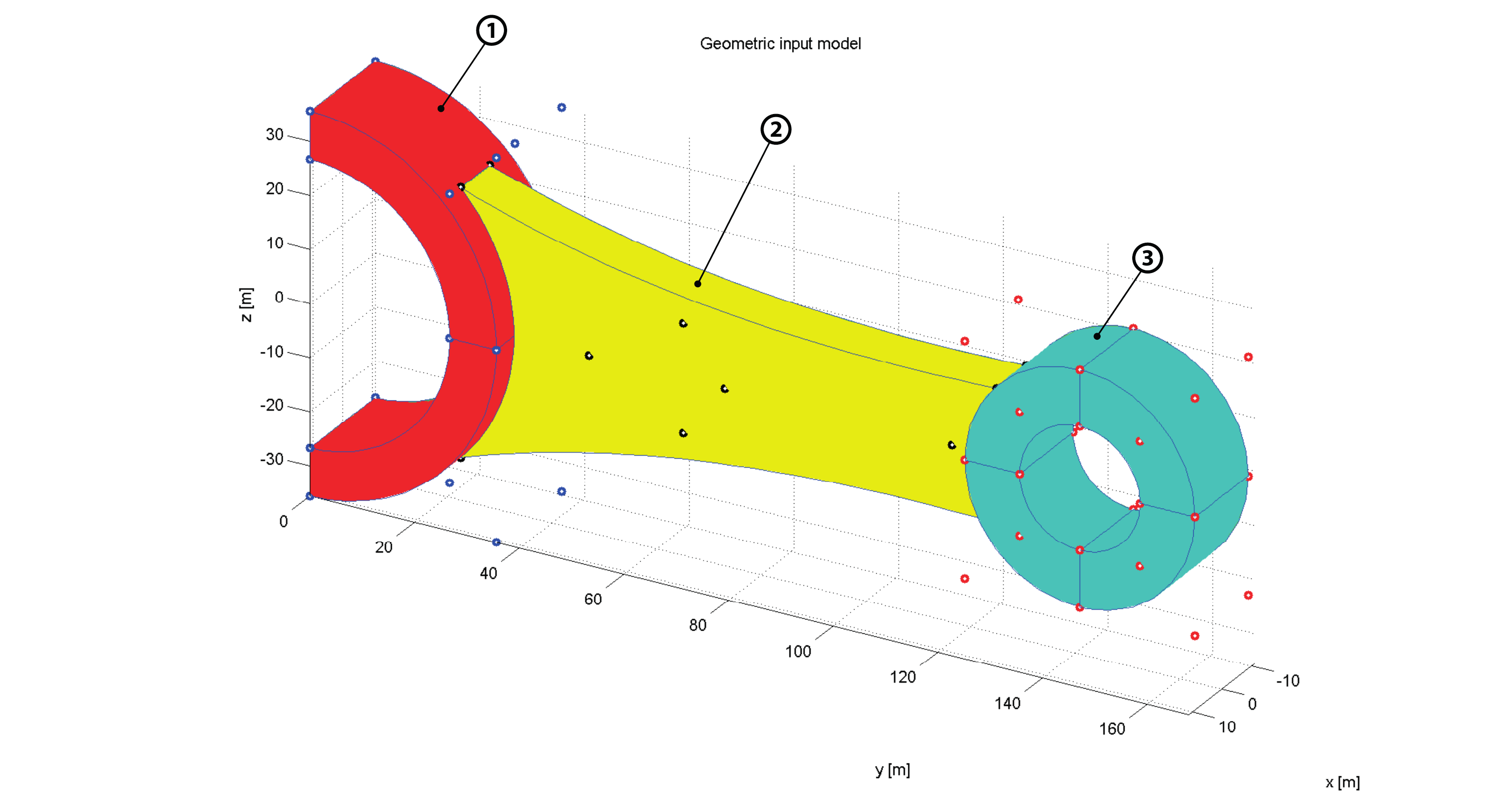}
  \caption{A multi-patch NURBS solid.}
  \label{fig:concepts}
\end{figure}

Nitsche's method \cite{nitsche} was originally proposed to weakly enforce Dirichlet boundary conditions
as an alternative to equivalent pointwise constraints.
The idea behind a Nitsche based approach is 
to replace the Lagrange multipliers arising in a dual formulation through their physical representation, 
namely the normal flux at the interface. Nitsche also added an extra penalty like term to restore 
the coercivity of the bilinear form. The method can be seen to lie in between the Lagrange multiplier
method and the penalty method.
The method has seen a resurgence in recent years and was applied for interface problems 
\cite{Hansbo20025537,Dolbow2009a},  for connecting overlapping meshes 
\cite{MZA:8203296,MZA:8203286,Sanders2012a,Sanders2011a}, for imposing Dirichlet boundary conditions 
in meshfree methods \cite{FernándezMéndez20041257}, in immersed boundary methods 
\cite{ruess2013weakly,NME:NME3339,embar_imposing_2010}, in fluid mechanics 
\cite{Bazilevs200712}, in the Finite Cell Method \cite{NME:NME4522}
and for contact mechanics \cite{nitsche-wriggers}. It has also been applied for stabilising constraints
in enriched finite elements \cite{Sanders2008a}.

The remainder of the paper is organised as follows. The problem description, governing equations
and weak formulation are presented in Section \ref{sec:problem}. Section \ref{sec:discretisation}
discusses the discretisation followed by implementation aspects given in Section \ref{sec:implementation}.
Several two and three dimensional examples are given in Section \ref{sec:examples}.

We denote $d_p$ and $d_s$ as the number of parametric directions and spatial directions respectively.
Both tensor and matrix notations are used.
In tensor notation, tensors of order one or greater are written in boldface. Lower case bold-face letters 
are used for first-order tensor whereas upper case bold-face letters indicate high-order tensors. 
The major exception to this rule are the physical second order stress tensor and the strain tensor
which are written in lower case. In matrix notation, the same symbols as for tensors are used to denote
the matrices but the connective operator symbols are skipped. 

\section{Problem description, governing equations and weak form}\label{sec:problem}

\subsection{Governing equations}

We define the domain $\Omega \subset \mathbb{R}^{d_s}$ with boundary $\Gamma \equiv \partial \Omega$. 
For sake of simplicity, we assume 
there is only one internal boundary denoted by $\Gamma_*$ that divides the domain into
two non-overlapping domains $\Omega^m, m=1,2$ such that $\Omega=\Omega^1 \cup \Omega^2$. 
In the context of multi-patch NURBS IGA, each domain represents a NURBS patch.
Excluding $\Gamma_*$, the rest of $\Gamma$ can be divided into Dirichlet and Neumann parts on each
domain, $\Gamma_u^m$ and $\Gamma_t^m$ respectively.
A superscript, $m$, is used to denote a quantity that is valid over region $\Omega^m$, with $m = 1, 2$.

With the primary unknown displacement field $\vm{u}^m$, the governing equations of linear elastostatic problems
are 
\begin{subequations}
\begin{alignat}{2}
-\nabla \;\bsym{\sigma}^m &=  \vm{b}^m         &\quad\text{on} \quad \Omega^m \\ 
             \vm{u}^m       &=  \bar{\vm{u}}^m   &\quad\text{on} \quad \Gamma_u^m \\
  \bsym{\sigma}^m \cdot \vm{n}^m &=  \bar{\vm{t}}^m  &\quad\text{on} \quad \Gamma_t^m \label{eq:Neumann}\\
  \vm{u}^1 &=  \vm{u}^2                        &\quad\text{on} \quad \Gamma_* \\
  \bsym{\sigma}^1 \cdot \vm{n}^1 &= -\bsym{\sigma}^2 \cdot \vm{n}^2  &\quad\text{on} \quad \Gamma_* \label{eq:tr}
\end{alignat}
\end{subequations}
where  $\bsym{\sigma}^m$ denotes the stress field; 
the last two equations express the continuity of displacements and tractions across $\Gamma_*$.
The prescribed displacement and traction are denoted by $\bar{\vm{u}}^m$ and $\bar{\vm{t}}^m$, respectively.
The outward unit normals to $\Omega^1$ and $\Omega^2$ are $\vm{n}^1$ and $\vm{n}^2$, respectively.

Under the small strain condition, the infinitesimal strain tensor reads 
$\bsym{\epsilon}^m=0.5(\nabla\vm{u}^m+\nabla\trans\vm{u}^m)$.
Constitutive equations are given by

\begin{equation}
\bsym{\sigma}^m = \vm{C}^m : \bsym{\epsilon}^m, \quad m=1,2
\end{equation}
where the constitutive tensors are denoted by $\vm{C}^1$ and $\vm{C}^2$. For linear 
isotropic elastic materials, the constitutive tensor is written as

\begin{equation}
C_{ijkl} = \lambda \delta_{ik}\delta_{kl} + \mu (\delta_{ik}\delta_{jl}+\delta_{il}\delta_{jk})
\end{equation}
where $\lambda=\frac{E\nu}{(1+\nu)(1-2\nu)}$ and $\mu=\frac{E}{2(1+\nu)}$ are the Lam{\'e} constants;
$E$ and $\nu$ are the Young's modulus and Poisson's ratio, respectively and $\delta_{ij}$ is the Kronecker
delta tensor. 

\begin{figure}
  \centering
  \includegraphics[width=\textwidth]{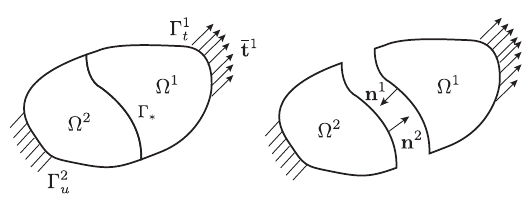}
  \caption{Computational domain with an internal interface.}
  \label{fig:domain}
\end{figure}

\subsection{Weak form}

We start by defining the spaces, $\bsym{S}^m$ and $\vm{V}^m$ over  domain $\Omega^m$ that will contain 
the solution and trial functions respectively:

\begin{equation}
\begin{split}
\bsym{S}^m&=\{\vm{u}^m(\vm{x})|\vm{u}^m(\vm{x}) \in \bsym{H}^1(\Omega^m), \vm{u}^m=\bar{\vm{u}}^m \;\;
\text{on $\Gamma_u^m$} \}\\
\bsym{V}^m&=\{\vm{w}^m(\vm{x})|\vm{w}^m(\vm{x}) \in \bsym{H}^1(\Omega^m), \vm{w}^m={\vm{0}} \;\;\text{on $\Gamma_u^m$} \}
\end{split}
\end{equation}

The standard application of Nitsche's method for the coupling is: 
Find $(\vm{u}^1,\vm{u}^2) \in \bsym{S}^1 \times \bsym{S}^2$ such that

\begin{multline}
\sum_{m=1}^2\int_{\Omega^m} \bsym{\epsilon}(\vm{w}^m):\bsym{\sigma}^m   \mathrm{d}\Omega 
-\int_{\Gamma^*}  \left(\jump{\vm{w}} \otimes \vm{n}^1\right) : \{\bsym{\sigma}\} \mathrm{d}\Gamma 
-\int_{\Gamma^*}  \left(\jump{\vm{u}} \otimes \vm{n}^1\right) : \{\bsym{\sigma}(\vm{w})\} \mathrm{d}\Gamma \\+
\int_{\Gamma^*}  \alpha \jump{\vm{w}} \cdot \jump{\vm{u}} \mathrm{d}\Gamma 
=  \sum_{m=1}^2 \int_{\Omega^m} \vm{w}^m\cdot\vm{b}^m \mathrm{d}\Omega + 
   \sum_{m=1}^2 \int_{\Gamma_t^m} \vm{w}^m \cdot \bar{\vm{t}}^m  \mathrm{d}\Gamma 
\label{eq:weakform}
\end{multline}
for all $(\vm{w}^1,\vm{w}^2) \in \bsym{V}^1 \times \bsym{V}^2$.
Derivation of this weak form is standard and can be found in, for example, \cite{Sanders2011a}. Note that we have assumed that
essential boundary conditions are enforced point-wise if possible or by other methods than Nitsche's method for we want to focus on the patch
coupling.

In Equation~\eqref{eq:weakform}, the jump and average operators, on the interface $\Gamma^*$, $\jump{\cdot}$
and $\{\cdot\}$ are defined as

\begin{equation}
\jump{\vm{u}} = \vm{u}^1 - \vm{u}^2, \quad
 \{\bsym{\sigma}\} = \frac{1}{2}(\bsym{\sigma}^1 + \bsym{\sigma}^2)
   \label{eq:jump-average}
\end{equation}
For completeness, note that the average operator for the stress field can be 
written generally as \cite{Sanders2012a}

\begin{equation}
\{\bsym{\sigma}\}=\gamma \bsym{\sigma}^1 + (1-\gamma)\bsym{\sigma}^2
\label{eq:general-average}
\end{equation}
where $0 \le \gamma \le 1$. The usual average operator is reproduced if $\gamma=0.5$ is used.
Equation~\eqref{eq:general-average} is often utilized to join a soft model and a stiff one \cite{Sanders2011a}.
Taking $\gamma=1$ (or $\gamma=0$) results in the one-sided mortaring method. 
In this paper, the standard average operator is used unless otherwise stated.

Except the second and third terms in the left hand side, Equation~\eqref{eq:weakform} is the same as the penalty
method. As in the penalty method, $\alpha$ is a free parameter for Nitsche's method. 
However, rather than being a penalty parameter, it should be viewed as a stabilization parameter 
in the context of this method. 
It has been shown \cite{Griebel} that a minimum $\alpha$ exists that will guarantee the positive definiteness 
of the bilinear form associated with Nitsche's method, thus, the stability of the method.

For discretisation we rewrite Equation~\eqref{eq:weakform} in a matrix form as follows:
Find $(\vm{u}^1,\vm{u}^2) \in \bsym{S}^1 \times \bsym{S}^2$ such that
\begin{multline}
 \sum_{m=1}^2\int_{\Omega^m}  (\bsym{\epsilon}(\vm{w}^m))\trans  \bsym{\sigma}^m  \mathrm{d}\Omega -
\int_{\Gamma_*}  \jump{\vm{w}}\trans \vm{n} \{\bsym{\sigma}\} \mathrm{d}\Gamma -
\int_{\Gamma_*}  \{\bsym{\sigma}(\vm{w})\}\trans \vm{n}\trans \jump{\vm{u}} \mathrm{d}\Gamma \\+
\int_{\Gamma_*}  \alpha \jump{\vm{w}}\trans \jump{\vm{u}} \mathrm{d}\Gamma 
=   \sum_{m=1}^2 \int_{\Gamma_t^m}(\vm{w}^m)\trans \bar{\vm{t}}^m  \mathrm{d}\Gamma 
+   \sum_{m=1}^2 \int_{\Omega^m}  (\vm{w}^m)\trans \vm{b}^m  \mathrm{d}\Omega
\label{eq:dg-weakform-matrix}
\end{multline}
for all $(\vm{w}^1,\vm{w}^2) \in \bsym{V}^1 \times \bsym{V}^2$.
Superscript T denotes the transpose operator.
Second order tensors ($\sigma_{ij}$ and $\epsilon_{ij}$) are written using the Voigt notation
as column vectors; $\bsym{\sigma}=[\sigma_{xx}, \sigma_{yy}, \sigma_{zz}, \sigma_{xy}, \sigma_{yz}, \sigma_{xz}]\trans$, $\bsym{\epsilon}=[\epsilon_{xx}, \epsilon_{yy}, \epsilon_{zz}, 2\epsilon_{xy}, 2\epsilon_{yz}, 2\epsilon_{xz}]\trans$,
and $\vm{n}$ (note that we removed the subscript 1 for subsequent derivations) 
is a matrix that reads

\begin{equation}
\vm{n}_{2D} = \begin{bmatrix}
n_x & 0 & n_y \\ 0 & n_y & n_x
\end{bmatrix}, \quad
\vm{n}_{3D} = \begin{bmatrix}
n_x & 0 & 0 & n_y & 0 & n_z \\
0 & n_y & 0 & n_x & n_z & 0\\
0 & 0 & n_z & 0 & n_y & n_x  
\end{bmatrix}\label{eq:n-matrix}
\end{equation}
for two dimensions and three dimensions, respectively.

\section{Discretisation}\label{sec:discretisation}

\subsection{NURBS}\label{sec:nurbs}

In this section, NURBS are briefly reviewed. 
We refer to the standard textbook \cite{piegl_book} for details.
A knot vector is a sequence in ascending order
of parameter values, written $\Xi=\{\xi_1,\xi_2,\ldots,\xi_{n+p+1}\}$
where $\xi_i$ is the \textit{i}th knot, $n$ is the number of basis functions and $p$ is 
the order of the B-spline basis. Open knots in which the first and last knots appear $p+1$ times are 
standard in the CAD literature and thus used in this manuscript \ie 
$\Xi=\{\underbrace{\xi_1,\ldots,\xi_1}_{\text{$p+1$ times}},\xi_2,\ldots,
\underbrace{\xi_m,\ldots\xi_m}_{\text{$p+1$ times}}\}$.

Given a knot vector $\Xi$, the B-spline basis functions are
defined recursively starting with the zeroth order basis
function ($p=0$) given by 

\begin{equation}
  N_{i,0}(\xi) = \begin{cases}
    1 & \textrm{if $ \xi_i \le \xi < \xi_{i+1}$}\\
    0 & \textrm{otherwise}
  \end{cases}
  \label{eq:basis-p0}
\end{equation}

\noindent and for a polynomial order $p \ge 1$

\begin{equation}
  N_{i,p}(\xi) = \frac{\xi-\xi_i}{\xi_{i+p}-\xi_i} N_{i,p-1}(\xi)
               + \frac{\xi_{i+p+1}-\xi}{\xi_{i+p+1}-\xi_{i+1}}
	       N_{i+1,p-1}(\xi)
  \label{eq:basis-p}
\end{equation}

\noindent This is referred to as the Cox-de Boor recursion formula.
Note that when evaluating these functions, ratios of the form $0/0$ are
defined as zero.

Some salient properties of B-spline basis functions are (1) they constitute 
a partition of unity, (2) each basis function is nonnegative over the entire domain,
(3) they are linearly independent, (4) the support of a B-spline function of order $p$ is $p+1$ knot spans
\ie $N_{i,p}$ is non-zero over $[\xi_i,\xi_{i+p+1}]$, (5) basis functions of order $p$ have $p-m_i$ continuous derivatives across knot $\xi_i$ where $m_i$ is the multiplicity of knot $\xi_i$ and (6) 
B-spline basis are generally only approximants (except at the ends of the parametric space interval, 
$[\xi_1,\xi_{n+p+1}]$) and not interpolants.

Fig.~\ref{fig:bspline-quad-open} illustrates a corresponding set of basis functions for an open, non-uniform knot vector. Of particular note is the interpolatory nature of the basis function at each
end of the interval created through an open knot vector, and the reduced continuity at $\xi = 4$ due to the presence of the location of a repeated knot where $C^0$ continuity is attained. Elsewhere, the functions are $C^1$ continuous ($C^{p-1}$).

\begin{figure}[h!]
  \centering 
  \includegraphics[width=0.7\textwidth]{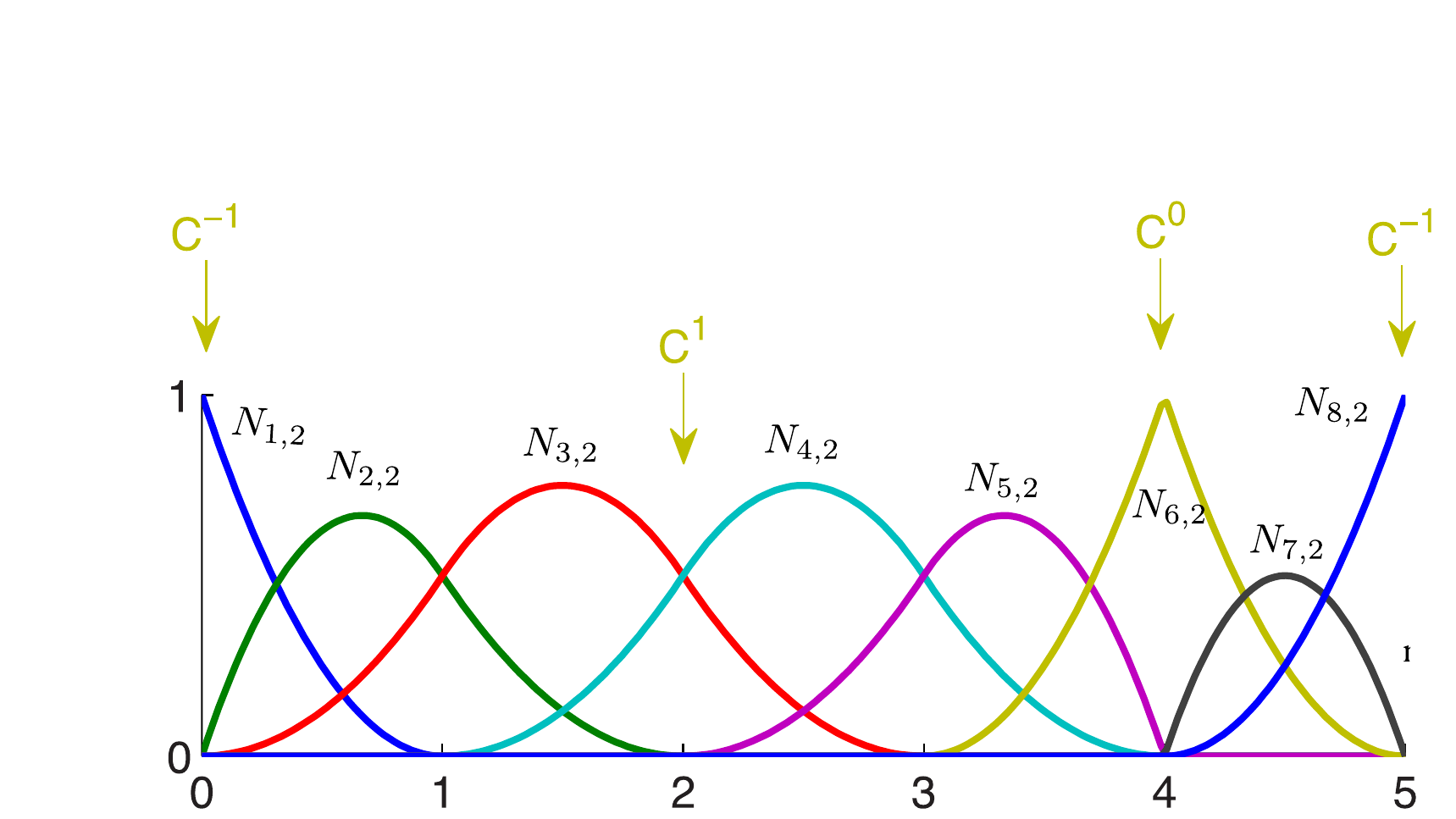}
  \caption{Quadratic B-spline basis functions defined for the open, non-uniform knot vector
  $\Xi=\{0,0,0,1,2,3,4,4,5,5,5\}$. Note the flexibility in the construction of
  basis functions with varying degrees of regularity.} 
  \label{fig:bspline-quad-open} 
\end{figure}

NURBS basis functions are defined as
\begin{equation}
	R_{i,p}(\xi) = \frac{N_{i,p}(\xi)w_i}{W(\xi)} =
	\frac{N_{i,p}(\xi)w_i}{\sum_{j=1}^{n}N_{j,p}(\xi)w_j}
  \label{eq:rational-basis}
\end{equation}
where $N_{i,p}(\xi)$ denotes the $i$th B-spline basis function of
order $p$ and $w_i$ are a set of $n$ positive weights. 
Selecting appropriate values for the $w_i$ permits the description of many
different types of curves including polynomials and circular arcs.
For the special case in which $w_i=c, i=1,2,\ldots,n$ the
NURBS basis reduces to the B-spline basis. Note that for simple geometries,
the weights can be defined analytically see \eg \cite{piegl_book}. For complex
geometries, they are obtained from CAD packages such as Rhino \cite{rhino}.

Let $\Xi^1=\{\xi_1,\xi_2,\ldots,\xi_{n+p+1}\}$, 
    $\Xi^2=\{\eta_1,\eta_2,\ldots,\eta_{m+q+1}\}$,
and $\Xi^3=\{\zeta_1,\zeta_2,\ldots,\zeta_{l+r+1}\}$ are the knot vectors
    and a control net $\vm{P}_{i,j,k} \in \mathds{R}^{d_s}$.
A tensor-product NURBS solid is defined as 

\begin{equation}
	\vm{V}(\xi,\eta,\zeta) = \sum_{i=1}^{n}\sum_{j=1}^{m}\sum_{k=1}^l
	\vm{P}_{i,j,k}  R_{i,j,k}^{p,q,r}(\xi,\eta,\zeta)
  \label{eq:NURBS-solid1}
\end{equation}
where the trivariate NURBS basis functions $R_{i,j,k}^{p,q,r}$ are given by
\begin{equation}
  R_{i,j,k}^{p,q,r}(\xi,\eta,\zeta) = \frac{N_{i}(\xi) M_j(\eta) P_k(\zeta) w_{i,j,k}}{
  \sum_{\hat{i}=1}^{n} \sum_{\hat{j}=1}^{m}  \sum_{\hat{k}=1}^{l}N_{\hat{i}}(\xi) M_{\hat{j}}(\eta)
  P_{\hat{k}}(\zeta) w_{\hat{i},\hat{j},\hat{k}}}.
  \label{eq:rational-basis3}
\end{equation}
By defining a global index $A$ through
\begin{equation}
  \label{eq:bspline_volume_mapping}
  A = (n \times m) ( k - 1)  + n( j - 1 ) + i
\end{equation}
a simplified form of Equation~\eqref{eq:NURBS-solid1} can be written as
\begin{equation}
  \label{eq:bspline_vol_simple}
  	\vm{V}(\boldsymbol{\xi}) = \sum_{A=1}^{n \times m \times l} \vm{P}_A  R_{A}^{p,q,r}(\boldsymbol{\xi} )
\end{equation}

\subsection{Isogeometric analysis}

\begin{figure}
  \centering
  \includegraphics[width=0.4\textwidth]{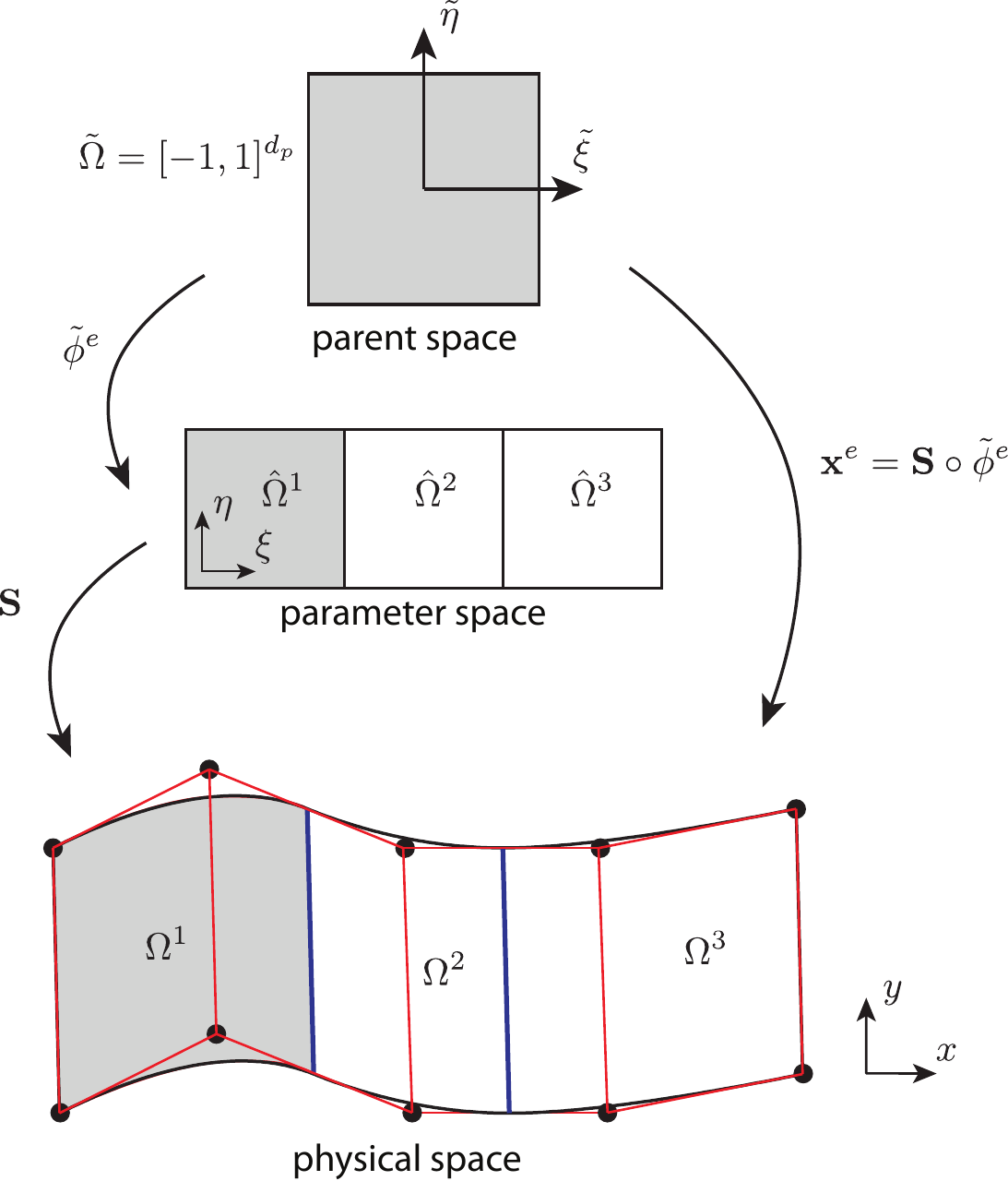}
  \caption{Diagrammatic interpretation of mappings from parent space ($\tilde{\Omega}$) 
     through parametric space ($\hat{\Omega}$) to physical space ($\Omega$). The parent space is where
     numerical quadrature rules are defined.}
  \label{fig:iga_mappings}
\end{figure}

Isogeometric analysis also makes use of an isoparametric formulation, but a key difference over its Lagrangian counterpart is the use of basis functions generated by CAD to discretise both the geometry and unknown fields. 
In IGA, regions bounded by knot lines with non-zero parametric area lead to a natural definition of 
element domains.
The use of NURBS basis functions for discretisation introduces the concept of parametric space which is absent in conventional FE implementations. The consequence of this additional space is that an additional mapping must be performed to operate in parent element coordinates. As shown in Fig.~\ref{fig:iga_mappings}, two mappings are considered for IGA with NURBS: a mapping $\tilde{\phi}^e: \tilde{\Omega} \to \hat{\Omega}^e$ and $\vm{S}: \hat{\Omega} \to \Omega$. The mapping $\vm{x}^e: \tilde{\Omega} \to \Omega^e$ is given by the composition $\vm{S}\circ \tilde{\phi}^e$. 

For a given element $e$, the geometry is expressed as
\begin{equation}
  \label{eq:iga_geometry_discretisation}
  \mathbf{x}^e(\tilde{\boldsymbol{\xi}}) = \sum_{a=1}^{n_{en}} \vm{P}_a^e R_a^e(\tilde{\boldsymbol{\xi}})
\end{equation}
where $a$ is a local basis function index, $n_{en} = (p+1)^{d_p}$ is the number of non-zero basis functions over element $e$ and $\vm{P}_a^e$,$R_a^e$ are the control point and NURBS basis function  associated with index $a$ respectively. We employ the commonly used notation of an element connectivity mapping \cite{hughes-fem-book} which translates a local basis function index to a global index through
\begin{equation}
  \label{eq:element_connectivity_array}
  A = \textrm{IEN}( a, e )
\end{equation}
Global and local control points are therefore related through $\vm{P}_A \equiv \vm{P}_{\textrm{IEN}(a,e)} \equiv \vm{P}_a^e$ with similar expressions for $R_a^e$.  

Taking the case $d_p = d_s = 2$,  an element defined by $\hat{\Omega}^e = [\xi_i, \xi_{i+1}]\otimes [\eta_i, \eta_{i+1}]$ is mapped from parent space to parametric space through 
\begin{align}
  \label{eq:phi_mapping}
  \tilde{\phi}^e(\tilde{\boldsymbol{\xi}})  =
\left\{ 
\begin{matrix}
\frac{1}{2}[(\xi_{i+1}-\xi_i)\tilde{\xi} + (\xi_{i+1}+\xi_i)]\\
 \frac{1}{2}[(\eta_{j+1}-\eta_j)\tilde{\eta} + (\eta_{j+1}+\eta_j)]
\end{matrix}
\right\}
\end{align}

A field $\vm{u}(\mathbf{x})$ which governs our relevant PDE can also be discretised in a similar manner to Equation~\eqref{eq:iga_geometry_discretisation} as
\begin{equation}
  \label{eq:iga_field_discretisation}
  \vm{u}^e(\tilde{\boldsymbol{\xi}}) =  \sum_{a=1}^{n_{en}} \vm{d}_a^e R_a^e(\tilde{\boldsymbol{\xi}})
\end{equation}
where $\vm{d}^e_a$ represents a control (nodal) variable. In contrast to conventional discretisations, these coefficients are not in general interpolatory at nodes. This is similar to the case of meshless
methods built on non-interpolatory shape functions such as the moving least squares
(MLS) \cite{efg-nayroles,NME:NME1620370205,nguyen_meshless_2008}. Using the Bubnov-Galerkin method, an analog
expansion as Equation~\eqref{eq:iga_field_discretisation} is adopted for the weight function and upon substituting
them into a weak form, a standard system of linear equations is obtained from which $\vm{d}$--the nodal variables
are obtained.

\subsection{Discrete equations}

The two domains $\Omega^m$ are discretised independently using finite elements.
At the interface $\Gamma_*$ there is a mismatch between the two meshes, cf. Fig.~\ref{fig:domain-mesh}.
The approximation of the displacement field is given by

\begin{figure}
  \centering
  \includegraphics[width=0.4\textwidth]{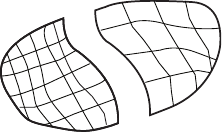}
  \caption{Independent discretisations of the domains.}
  \label{fig:domain-mesh}
\end{figure}

\begin{equation}
  \vm{u}^m  =  N_A^m \vm{a}_A^m
  \label{eq:approximation}
\end{equation}
where $N^m_A$ denotes the finite element shape functions associated to domain $\Omega^m$ (which can be 
      any Lagrange shape functions or the B-spline and NURBS basis functions presented in Section \ref{sec:nurbs}) and
$\vm{a}_A^m=[a_{xA}^m\; a_{yA}^m\; a_{zA}^m]\trans$ represents the nodal displacements of domain $\Omega^m$.

The stresses, strains and displacements are given by 

\begin{equation}
\bsym{\sigma}^m = \vm{C}^m \vm{B}^m \vm{a}^m, \quad
\bsym{\epsilon}^m = \vm{B}^m \vm{a}^m, \quad 
\vm{u}^m=\vm{N}^m \vm{a}^m
\label{eq:11}
\end{equation}
where $\vm{B}$ is the standard strain-displacement matrix and $\vm{N}$ represents the standard
shape function matrix. For two dimensional element $e$, they are given by

\begin{equation}
\vm{B}_e^m = \begin{bmatrix}
N_{1,x}^m & 0 & N_{2,x}^m & 0 & \ldots\\
0 & N_{1,y}^m & 0 & N_{2,y}^m & \ldots \\
N_{1,y}^m & N_{1,x}^m & N_{2,y}^m & N_{2,x}^m & \ldots 
\end{bmatrix},\quad \vm{N}_e^m = \begin{bmatrix}
N_1^m & 0 & N_2^m & 0 & \ldots\\
0 & N_1^m & 0 & N_2^m & \ldots       
\end{bmatrix}
\end{equation}
Expressions for three dimensional elements can be found in many FEM textbooks \eg
\cite{hughes-fem-book}. The notation $N_{I,x}$ denotes
the derivative of shape function $N_I$ with respect to $x$. This notation for partial derivatives
will be used in subsequent sections.

The jump operator and the average operator are given by

\begin{equation}
\begin{split}
\jump{\vm{u}}     &= \vm{N}^1\vm{a}^1 - \vm{N}^2 \vm{a}^2\\
\{\bsym{\sigma}\} &= \frac{1}{2}\left( \vm{C}^1 \vm{B}^1 \vm{a}^1 + \vm{C}^2\vm{B}^2 \vm{a}^2 \right)
\end{split}
\label{eq:12}
\end{equation}
and analog expansions are used for $\jump{\vm{w}}$ and  $\{\bsym{\sigma}(\vm{w})\} $
\begin{equation}
\begin{split}
\jump{\vm{w}}     &= \vm{N}^1 \delta \vm{a}^1 - \vm{N}^2 \delta \vm{a}^2\\
\{\bsym{\sigma}(\vm{w})\} &= \frac{1}{2}\left( \vm{C}^1 \vm{B}^1 \delta \vm{a}^1 + 
      \vm{C}^2\vm{B}^2 \delta \vm{a}^2 \right)
\end{split}
\label{eq:13}
\end{equation}

Upon substituting Equations~\eqref{eq:11},\eqref{eq:12} and \eqref{eq:13} into 
Equation~\eqref{eq:dg-weakform-matrix} and invoking the arbitrariness of $\delta \vm{a}^m$, we obtain
the discrete equation that can be written as

\begin{equation}
\left[\vm{K}^b + \vm{K}^n + (\vm{K}^n)\trans + \vm{K}^s\right] \vm{a} = \vm{f}_\text{ext}
\end{equation}
in which $\vm{K}^b$ denotes the bulk stiffness matrix; $\vm{K}^n$ and $\vm{K}^s$ are the interfacial
stiffness matrices or the coupling matrices. The external force vector is denoted by $\vm{f}_\text{ext}$
and its expression is standard and thus presented here.

The bulk stiffness matrix is given by

\begin{equation}
 \vm{K}^b = \sum_m^2 \int_{\Omega^m} (\vm{B}^m)\trans \vm{C}^m \vm{B}^m\di \Omega 
\end{equation}
and the coupling matrices are given by
\begin{equation}
\vm{K}^n = \begin{bmatrix}
-\D\int_{\Gamma_*} \vm{N}^{1\text{T}} \vm{n} \frac{1}{2}\vm{C}^1\vm{B}^1 \mathrm{d}\Gamma &
-\D\int_{\Gamma_*} \vm{N}^{1\text{T}} \vm{n} \frac{1}{2}\vm{C}^2\vm{B}^2 \mathrm{d}\Gamma \\
\D\int_{\Gamma_*} \vm{N}^{2\text{T}} \vm{n} \frac{1}{2}\vm{C}^1\vm{B}^1 \mathrm{d}\Gamma &
\D\int_{\Gamma_*} \vm{N}^{2\text{T}} \vm{n} \frac{1}{2}\vm{C}^2\vm{B}^2 \mathrm{d}\Gamma 
\end{bmatrix}\label{eq:nitsche-kdg}
\end{equation}
and by
\begin{equation}
\vm{K}^s = \begin{bmatrix}
\D\int_{\Gamma_*}  \alpha  \vm{N}^{1\text{T}} \vm{N}^1 \mathrm{d}\Gamma &
 - \D\int_{\Gamma_*}  \alpha  \vm{N}^{1\text{T}} \vm{N}^2 \mathrm{d}\Gamma\\
 - \D\int_{\Gamma_*}  \alpha  \vm{N}^{2\text{T}} \vm{N}^1 \mathrm{d}\Gamma&
 \D\int_{\Gamma_*}  \alpha  \vm{N}^{2\text{T}} \vm{N}^2 \mathrm{d}\Gamma
\end{bmatrix}\label{eq:nitsche-kpe}
\end{equation}
If the average operator defined in Equation~\eqref{eq:general-average} is used, 
we have   

\begin{equation}
\vm{K}^n = \begin{bmatrix}
-\gamma\D\int_{\Gamma_*} \vm{N}^{1\text{T}} \vm{n} \vm{C}^1\vm{B}^1 \mathrm{d}\Gamma &
-(1-\gamma)\D\int_{\Gamma_*} \vm{N}^{1\text{T}} \vm{n} \vm{C}^2\vm{B}^2 \mathrm{d}\Gamma \\
\gamma\D\int_{\Gamma_*} \vm{N}^{2\text{T}} \vm{n} \vm{C}^1\vm{B}^1 \mathrm{d}\Gamma &
(1-\gamma)\D\int_{\Gamma_*} \vm{N}^{2\text{T}} \vm{n} \vm{C}^2\vm{B}^2 \mathrm{d}\Gamma 
\end{bmatrix}\label{eq:nitsche-kdg-general}
\end{equation}

\section{Implementation}\label{sec:implementation}

For the computation of the bulk stiffness matrices is standard, in this section we focus on
the implementation of the coupling matrices for both two and three dimensional problems. 
For sake of presentation, Lagrange finite elements are discussed firstly and generalisation to
NURBS elements is given subsequently with minor modifications.

\subsection{Two dimensions}

\begin{figure}
  \centering
  \includegraphics[width=0.8\textwidth]{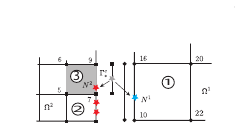}
  \caption{Independent discretisations of the domains: hierarchical meshes. The interface $\Gamma_*$ is
  discretised using the element edges of $\Omega^2$ that intersect $\Gamma_*$. For the grey element, the Gauss point is denoted by the red star which is mapped to the GP in element 1 (green star). }
  \label{fig:domain-mesh-hier}
\end{figure}

\subsubsection{Hierarchical meshes}

First, we consider hierarchical meshes as shown in Fig.~\ref{fig:domain-mesh-hier}.
In this case, the interface integrals can be straightforwardly calculated as explained in what follows.
Let assume that a fine mesh is adopted for $\Omega^2$ and a coarse mesh for $\Omega^1$, cf.
Fig.~\ref{fig:domain-mesh-hier}. We use the fine elements on $\Gamma_*$ to evaluate the interfacial integral

\begin{equation}
\int_{\Gamma_*} f(N^1,N^2) d \Gamma = \bigcup_{e=1}^{nbe} 
\int_{\Gamma_*^e} f(N^1,N^2) d \Gamma 
\end{equation}
where $\Gamma_*^e=\Gamma_* \cap \Omega_e^{2,b}$ and $\{\Omega_e^{2,b}\}_{1}^{nbe}$ denotes elements in $\Omega^2$
that intersect with $\Gamma_*$.
What makes hierarchical meshes attractive is that for a fine element on $\Gamma_*$ one knows
the element in the coarse mesh that locates the other side of the interface.

For the elemental interface integral, a Gauss quadrature rule for line 
elements is adopted. For example, two GPs are used for bilinear elements. Let the
GPs denoted by  $\{\xi_i\}_{i=1}^{ngp}$. These GPs have to be mapped to two parent elements--
one associated with $\Omega_e^{2,b}$ and  one associated with $\Omega_e^{1,b}$.
That is given $\xi_i$, one has to solve for $\bsym{\xi}_i^2$ and $\bsym{\xi}_i^1$
($\bsym{\xi}_i^2=(\xi_i^2,\eta_i^2)$)

\begin{equation}
\begin{split}
\vm{x}_i &= \vm{M}(\xi_i)\vm{x}_l \\
\vm{x}_i &= \vm{N}^2(\bsym{\xi}_i^2)\vm{x}_e^2 \rightarrow \bsym{\xi}_i^2\\
\vm{x}_i &= \vm{N}^1(\bsym{\xi}_i^1)\vm{x}_e^1 \rightarrow \bsym{\xi}_i^1
\end{split}
\end{equation}
where the first equation is used to compute the global coordinates of the GP ($\vm{x}_i=(x_i,y_i)$)
and the second and third equations are used to compute the natural coordinates
of the GP in the parent element associated with $\Omega_e^{k,b}$. 
Usually a Newton-Raphson method is used for this.
In the above,
$\vm{M}$ denotes the row vector of shape functions of a two-noded line element; $\vm{x}_l$
are the nodal coordinates of two boundary nodes of $\Gamma_*^e$ (for the example given in 
Fig.~\ref{fig:domain-mesh-hier}, they are nodes 7 and 9);
$\vm{x}_e^k$ ($k=1,2$) denotes the nodal coordinates of $\Omega_e^{k,b}$. 
$\vm{N}^k$ denote the row vector of shape functions of element $\Omega_e^{k,b}$.
For the example given in Fig.~\ref{fig:domain-mesh-hier}, $\vm{x}_e^2$ stores the
coordinates of nodes 5,7,9 and 6. And, $\vm{x}_e^1$ stores the coordinates of nodes 10,22,20 and 16.

It is now ready to evaluate the interfacial integral as

\begin{equation}
\int_{\Gamma_*^e} f(N^1,N^2) \di \Gamma = \sum_{i=1}^{ngp} f(N^1(\bsym{\xi}_i^1),N^2(
\bsym{\xi}_i^2) ) w_i
\end{equation}
where $w_i$ equals the weight multiplied with the Jacobian of the transformation from
the line parent element $[-1,1]$ to $\Gamma_*^e$.

Finally the coupling terms are assembled to the global stiffness matrix in a standard manner.
For example $\vm{K}^{n,11}$ is assembled using the connectivity of $\Omega_e^{1,b}$
and $\vm{K}^{n,22}$ is assembled using the connectivity of $\Omega_e^{2,b}$.

\subsubsection{Non-matching structured meshes}

Non-matching structured meshes are plotted in Fig.~\ref{fig:domain-mesh5}.
In those cases, the evaluation of the interfacial integrals are more complicated.
We use the trace mesh of $\Omega^1$ on the coupling interface $\Gamma_*$ to perform the
numerical integration. We use two data structures to store the Gauss points namely (for the
concrete example shown in Fig.~\ref{fig:domain-mesh5})      
$gp1=\{(\bsym{\xi}^1_i,w_i,e^1_i)\}_{i=1}^4$ and $gp2=\{(\bsym{\xi}^2_i,e^2_i)\}_{i=1}^4$ where
$e^m_i$ indicates the index of element of $\Omega^m$ that contains GP $\bsym{\xi}^m_i$.
After having these GPs, the assembly of the coupling matrices follows the procedure outlined in
Box \ref{box-k-example}.

  \begin{figure}
  \centering
  \includegraphics[width=0.45\textwidth]{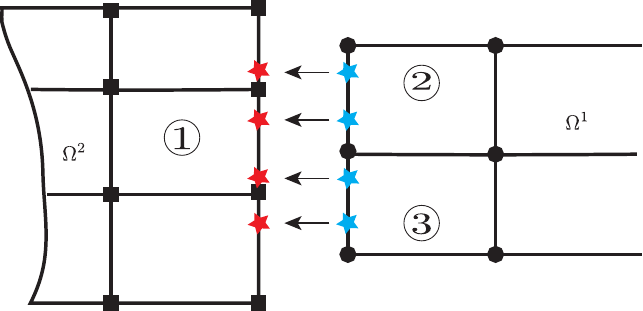}
  \caption{Independent discretisations of the domains: non-matching structured meshes.}
  \label{fig:domain-mesh5}
\end{figure}

\begin{Fbox}
	\caption{Assembly of coupling matrices}
  \begin{enumerate}
	  \item Loop over Gauss points (GPs), $i$ 
  \begin{enumerate}
	  \item Get $\bsym{\xi}^1_i$, $w_i$ and $e^1_i$ from $gp1$
	  \item Get $\bsym{\xi}^2_i$ and $e^2_i$ from $gp2$
	  \item Compute shape functions $\vm{N}^1(\bsym{\xi}^1_i)$
	  \item Compute shape functions $\vm{N}^2(\bsym{\xi}^2_i)$
	  \item Compute $\vm{K}^{s,12}= -\alpha  \vm{N}^{1\text{T}} \vm{N}^2 w_i$
	  \item Assemble  $\vm{K}^{s,12}$ to the global stiffness matrix using the connectivity array
          of $e^1_i$ (rows) and $e^2_i$ (columns).
  \end{enumerate}
  \item End loop over GPs
  \end{enumerate}
  \label{box-k-example}
\end{Fbox}

\subsection{Three dimensional formulations}

\begin{figure}
  \centering
  \includegraphics[width=0.95\textwidth]{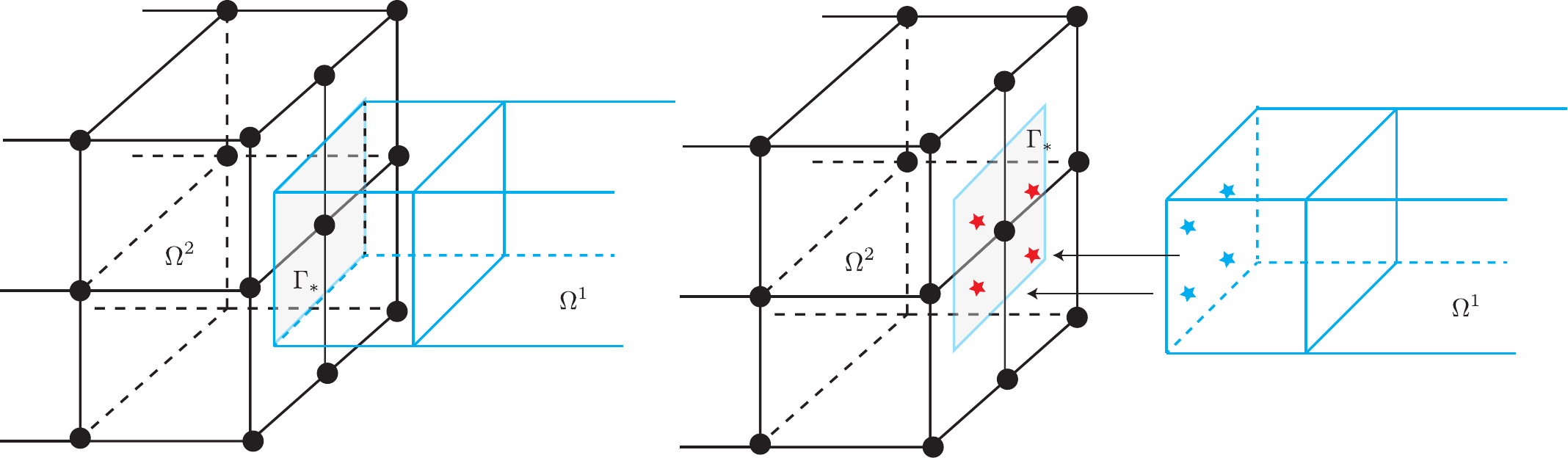}
  \caption{Coupling of two three dimensional continuum models. For evaluating the coupling terms,
  we use the trace mesh of $\Omega^1$ on the coupling interface $\Gamma_*$. In this figure, there is
  only one element of the trace mesh for sake of illustration.}
  \label{fig:3d-3d-coupling}
\end{figure}

This section presents the implementation for 3D, we refer to Fig.~\ref{fig:3d-3d-coupling}.
The computation of GPs required for the coupling matrices is given in Box \ref{box-gp-3D}.
After having obtained $gp1$ and $gp2$ data structures, the assembly of the coupling matrices
follows Box \ref{box-k-example}.

\begin{Fbox}[thpb]
	\caption{Determination of $gp1$ and $gp2$}
  \begin{enumerate}
	  \item For each element $e^1$ of the trace mesh, do  
  \begin{enumerate}
	  \item Distribute GPs on the face, $\{(\xi_i,\eta_i,w_i)\}_{i=1}^{ngp}$
	  \item Loop over the GPs, $i$
  \begin{enumerate}
	  \item Transform GP $i$ to physical space using
               \begin{equation}
               \label{eq:hang}
               \vm{x}_i = \vm{M}(\xi_i,\eta_i)\vm{x}_l 
               \end{equation}
           \item Compute tangent vectors, normal vector and the weight    
               \begin{equation}
               \vm{a}_1 = \vm{M}_{,\xi} \vm{x}_l,\quad
               \vm{a}_2 = \vm{M}_{,\eta} \vm{x}_l, 
               \quad \vm{n}=\frac{\vm{a}_1\times\vm{a}_2}{\norm{\vm{a}_1\times\vm{a}_2}}, \quad
               \bar{w}_i = w_i \norm{\vm{a}_1\times\vm{a}_2}
               \label{eq:tangents}
               \end{equation}
	  \item Transform GP $i$ from physical space to parent space of $\Omega^1$ using
                \begin{equation}
                \vm{x}_i = \vm{N}^1(\xi_i^1,\eta_i^1,\zeta_i^1)\vm{x}_e^1
                \rightarrow (\xi_i^1,\eta_i^1,\zeta_i^1) \label{eq:hang1}
                \end{equation}
          \item Find index of element in $\Omega^2$ that contains $\vm{x}_i$, named it $e^2$
	  \item Transform GP $i$ from physical space to parent space of $\Omega^2$ using
               \begin{equation}
               \vm{x}_i = \vm{N}^2 (\xi_i^2,\eta_i^2,\zeta_i^2) \vm{x}_e^2 \rightarrow 
               (\xi_i^2,\eta_i^2,\zeta_i^2) \label{eq:hang2}
               \end{equation}
               where $\vm{x}_e^2$ are the nodal coordinates of element $e^2$.
  \end{enumerate}
  \item End loop over GPs
  \end{enumerate}
  \item End for
  \end{enumerate}
  \label{box-gp-3D}
\end{Fbox}

\subsection{Extension to NURBS elements}
 
Since NURBS basis functions are defined on the parameter space not on the parent space, there is a slight modification to the implementation. The GPs are now give by  $\{(\tilde{\xi}_i,\tilde{\eta}_i,\tilde{w}_i)\}_{i=1}^{ngp}$.
They are firstly transformed to the parameter space using the mapping defined in Equation~\eqref{eq:phi_mapping}:
$\{(\xi_i,\eta_i,w_i)\}_{i=1}^{ngp}$ where $w_i = \tilde{w}_i J$ with $J$ is the Jacobian of the parent-to-parameter mapping. After that one works with the parameter space, for example the inverse mapping 
Equation~\eqref{eq:hang1} determines a point in the parameter space.

Steps (iv) and (v) in the algorithm given in Box \ref{box-gp-3D} demand modifications because one can exploit
the fact that the NURBS mapping, Equation~\eqref{eq:bspline_vol_simple}, is global. Hence, one writes Equation~\eqref{eq:hang2}
as follows               
               \begin{equation}
               \vm{x}_i = \vm{N}^2 (\xi_i^2,\eta_i^2,\zeta_i^2) \vm{x}^2 \rightarrow 
               (\xi_i^2,\eta_i^2,\zeta_i^2) \label{eq:hang3}
               \end{equation}
               where $\vm{x}^2$ are the control point of patch 2. Note that in Equation~\eqref{eq:hang1}, $\vm{x}^1_e$ denotes the
               control points of only the element under consideration. Using the output $ (\xi_i^2,\eta_i^2,\zeta_i^2)$ and the standard $FindSpan$ algorithm, cf. \cite{piegl_book},
               one can determine which element $\vm{x}_i$ belongs to \ie $e^2$.

\begin{rmk}
Note also that if B{\'e}zier extraction is used to implement NURBS-based
IGA, see \eg \cite{borden_isogeometric_2011}, then this section can be ignored since with B{\'e}zier extraction
the basis are the Bernstein basis, which are defined in the parent space as well, multiplied 
with some sparse matrices. Moreover, B{\'e}zier extraction will facilitate the incorporation of the non-conforming multi-patch
NURBS IGA into existing FE codes including commercially available FE packages.
\end{rmk}

\section{Numerical examples}\label{sec:examples}

In this section three numerical examples of increasing complexity are presented to assess the performance
of the proposed method. They are listed as follows

\begin{enumerate}
   \item Timoshenko beam (2D/2D coupling)
   \item Cantilever beam (3D/3D coupling)
   \item Connecting rod (complex 3D/3D coupling)
\end{enumerate}

The first two examples are simple problems to verify the implementation and we provide convergence analysis
for the first example.
Unless otherwise stated, we use 
MIGFEM--an open source Matlab IGA code which is available at \url{https://sourceforge.net/projects/cmcodes/}
for our computations and  the  visualisation was performed in Paraview \cite{para}.

\subsection{Timoshenko beam}

Consider a beam of dimensions $L \times D$ (unit thickness), subjected to a
parabolic traction at the free end as shown in Fig.~\ref{fig:beam-geo}.
A plane stress state is assumed. The parabolic traction is given by

\begin{equation}
    t_y(y) = -\frac{P}{2I} \biggl ( \frac{D^2}{4} - y^2 \biggr)\label{eq:ty}
\end{equation}

\noindent where $I = D^3/12$ is the moment of inertia. The exact displacement
field of this problem is, see \eg \cite{elasticity_book}

\begin{equation}
\begin{split}
    u_x(x,y) &=  \frac{Py}{6EI} \biggl [ (6L-3x)x + (2+\nu)\biggl(y^2-\frac{D^2}{4}\biggr) \biggr] \\
    u_y(x,y) &=  - \frac{P}{6EI} \biggl [ 3\nu y^2(L-x) + (4+5\nu)\frac{D^2x}{4} +(3L-x)x^2 \biggr] \\
\end{split}
\label{eq:tBeamExactDisp}
\end{equation}

\noindent and the exact stresses are

\begin{equation}
	\sigma_{xx}(x,y) =  \frac{P(L-x)y}{I}; 
	\quad \sigma_{yy}(x,y) =  0, \quad
\sigma_{xy}(x,y) =  -\frac{P}{2I} \biggl ( \frac{D^2}{4}-y^2\biggr)
\end{equation}

\noindent In the computations, material properties are taken as $E=
3.0 \times 10^7$, $\nu = 0.3$ and the beam dimensions are $D=6$ and
$L=48$. The shear force is $P = 1000$.  In order to model the clamping condition, 
the displacement defined by Equation~\eqref{eq:tBeamExactDisp} is prescribed as essential boundary 
conditions at $x=0, -D/2 \le y \le D/2$. This problem is solved with bilinear Lagrange elements
and high order B-splines elements. The former helps to verify the implementation in addition to
the ease of enforcement of Dirichlet boundary conditions (BCs). For the latter, care must be taken
in enforcing the Dirichlet BCs given in Equation~\eqref{eq:tBeamExactDisp} since the B-splines
are not interpolatory. The beam is divided into two domains by a vertical line at $x=L/2$ \ie
$\Gamma^*=\{x=L/2,-D/2 \le y \le D/2\}$. \\

\begin{figure}[htbp]
  \centering 
  \psfrag{p}{P}\psfrag{l}{$L$}\psfrag{d}{$D$}\psfrag{x}{$x$}\psfrag{y}{$y$}
  \includegraphics[width=0.5\textwidth]{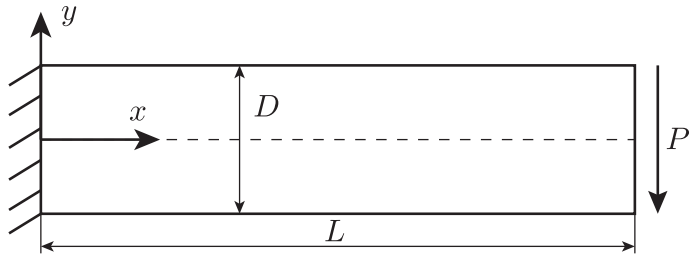}
  \caption{Timoshenko beam: problem description.}
  \label{fig:beam-geo}
\end{figure}

\noindent \textbf{Lagrange elements} Firstly, a conforming mesh (however there are double nodes at $\Gamma^*$) 
is considered and each domain is
discretised by a mesh of $20\times4$ elements as given in Fig.~\ref{fig:beam-meshes}a. Then, a non-conforming
mesh where the left domain is discretised by $20\times8$ elements and
the right domain is meshed by $20\times4$ is considered, cf. Fig.~\ref{fig:beam-meshes}b. 
A value of $1\times10^8$ was used for $\alpha$.
The vertical displacements along the midline of the beam 
($u_y(0 \le x \le L,y=0$) are plotted in Fig.~\ref{fig:beam-disp} together with the exact solution.
A good agreement can be observed. 
The stresses are plotted in Fig.~\ref{fig:beam-stress}.\\

\begin{figure}[htbp]
  \centering 
   \subfloat[Conforming mesh]{\includegraphics[width=0.45\textwidth]{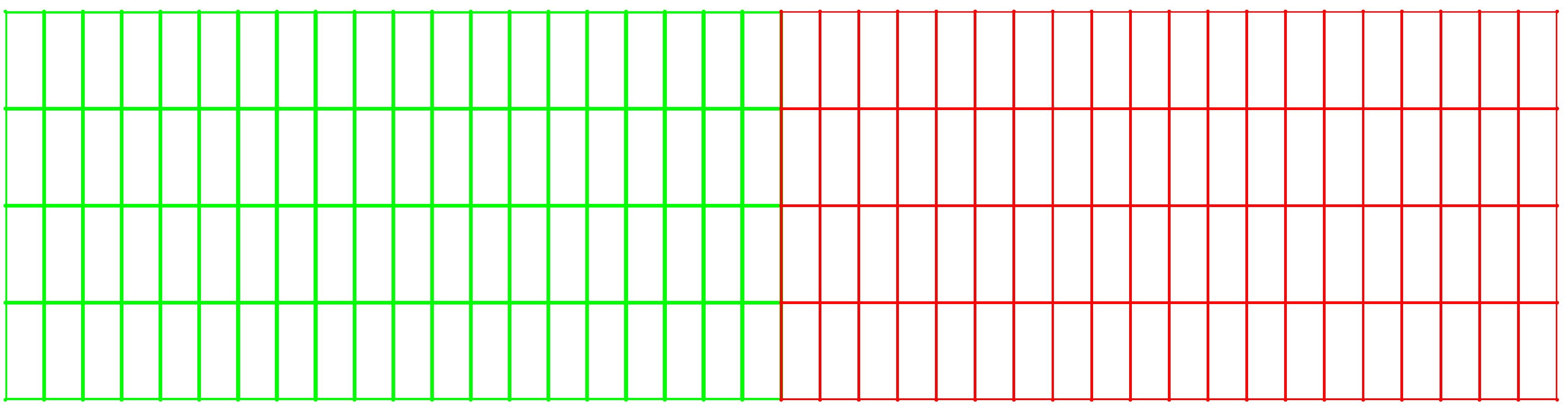}}\\
   \subfloat[Non conforming mesh]{\includegraphics[width=0.45\textwidth]{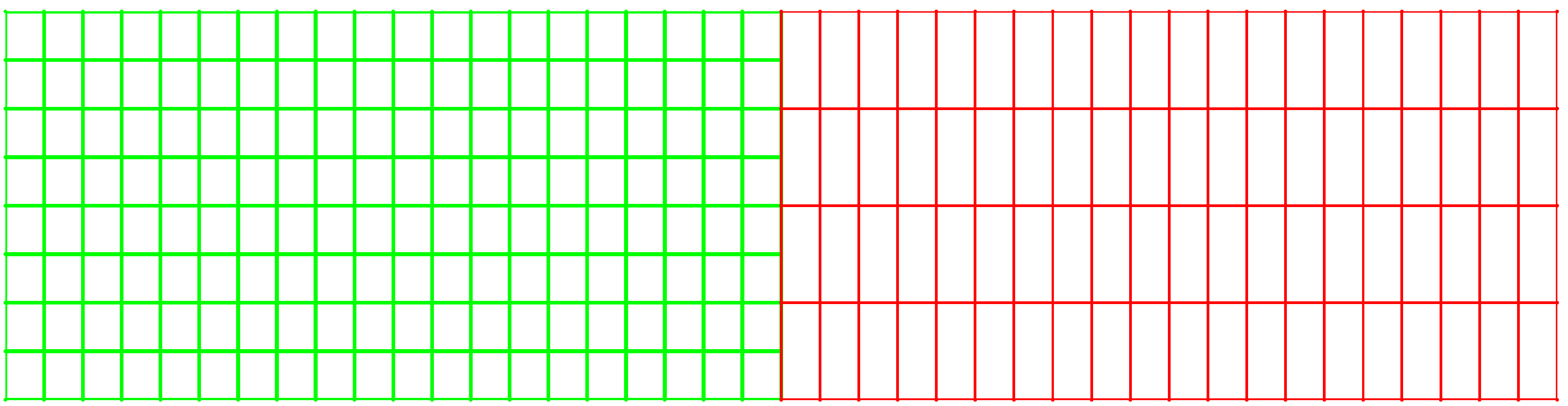}}
  \caption{Timoshenko beam: conforming and non-conforming meshes. Note that even with the conforming mesh,
  there are double nodes at the coupling interface $x=L/2, -D/2 \le y \le D/2$.}
  \label{fig:beam-meshes}
\end{figure}

\begin{figure}[htbp]
  \centering 
   \subfloat[Conforming mesh]{\includegraphics[width=0.45\textwidth]{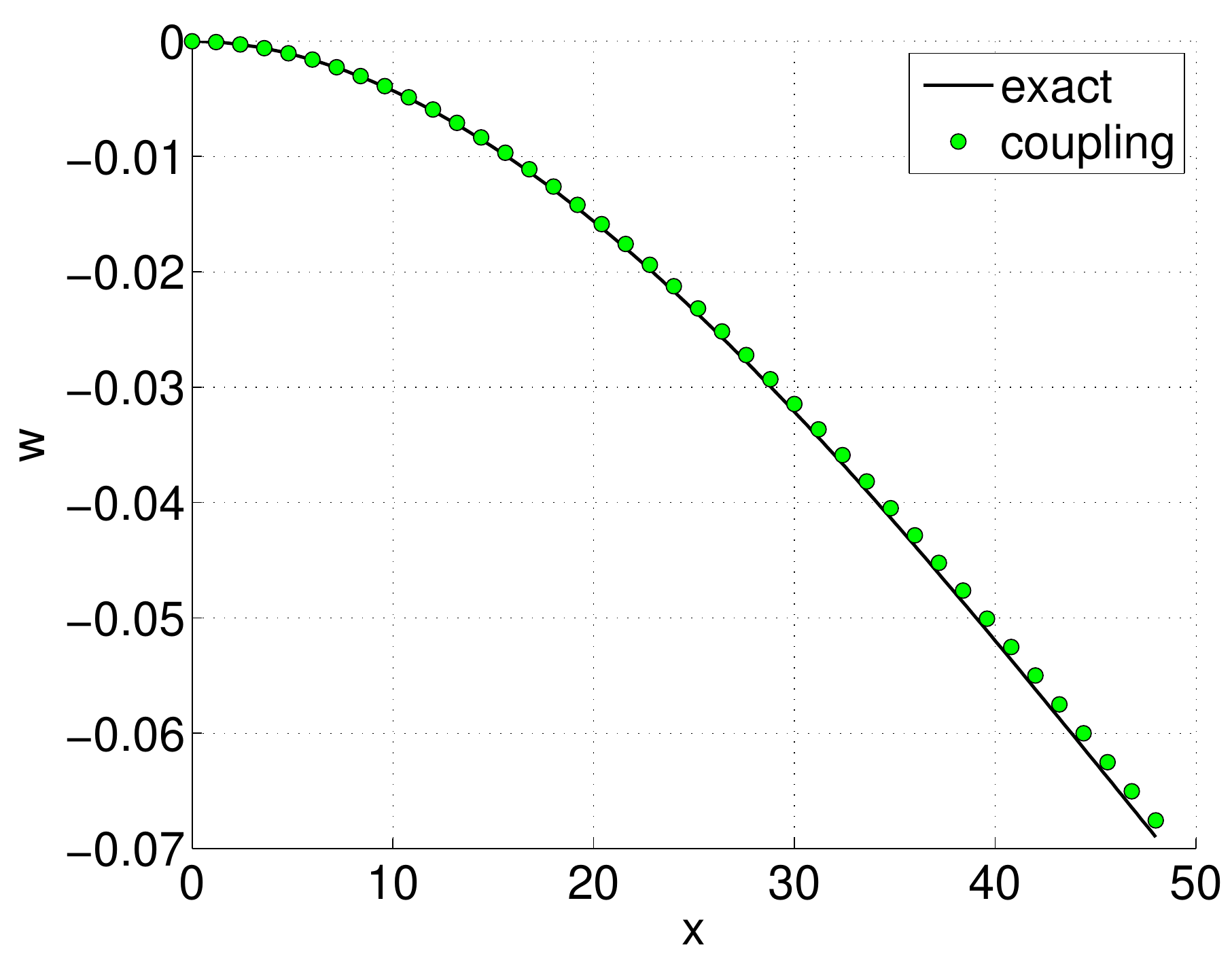}}
   \subfloat[Non-conforming mesh]{\includegraphics[width=0.45\textwidth]{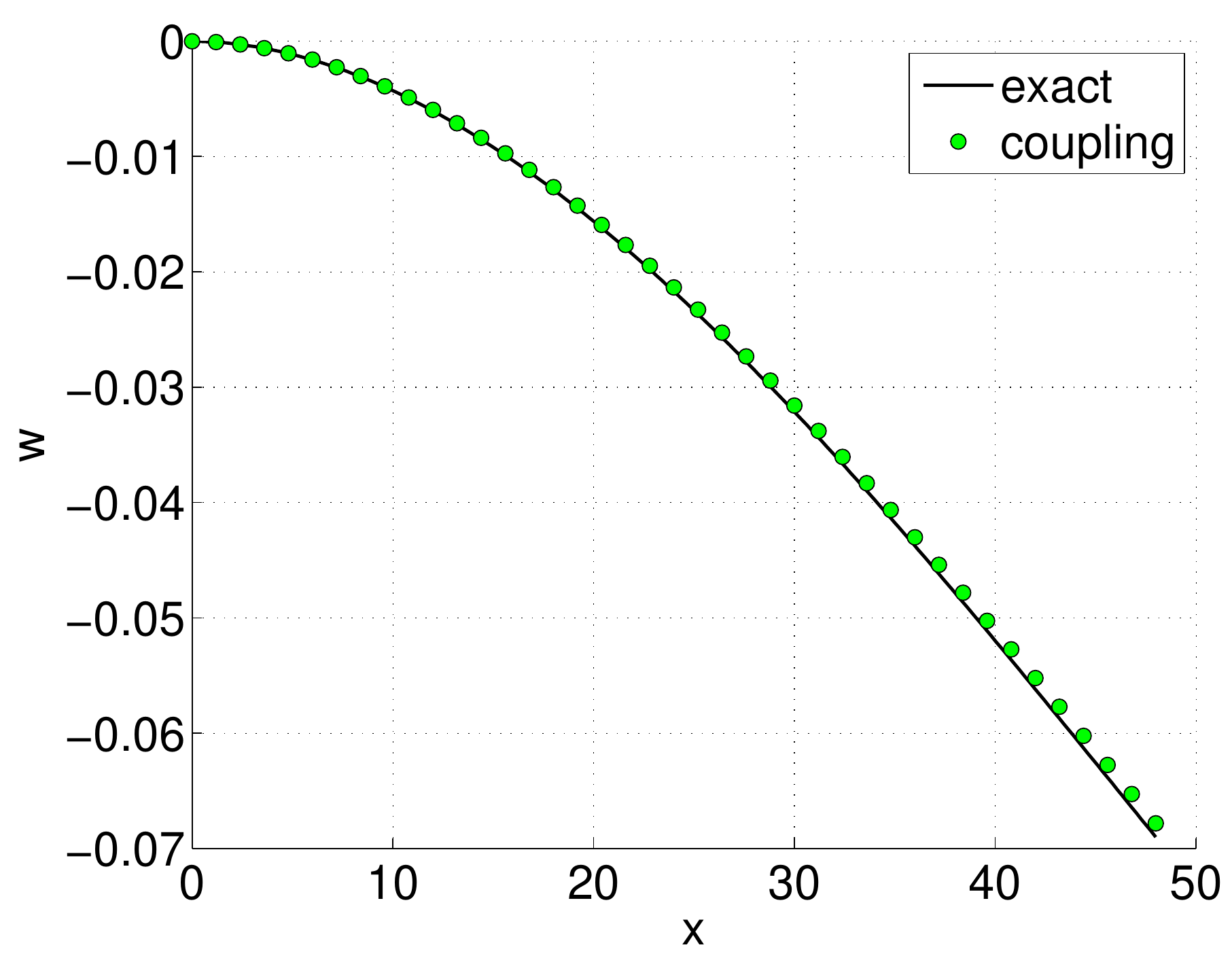}}
  \caption{Timoshenko beam: comparison of $u_y(0 \le x \le L,y=0)$ with the exact solution.}
  \label{fig:beam-disp}
\end{figure}

\begin{figure}[htbp]
  \centering 
   \subfloat[Stresses along the beam length]{\includegraphics[width=0.45\textwidth]{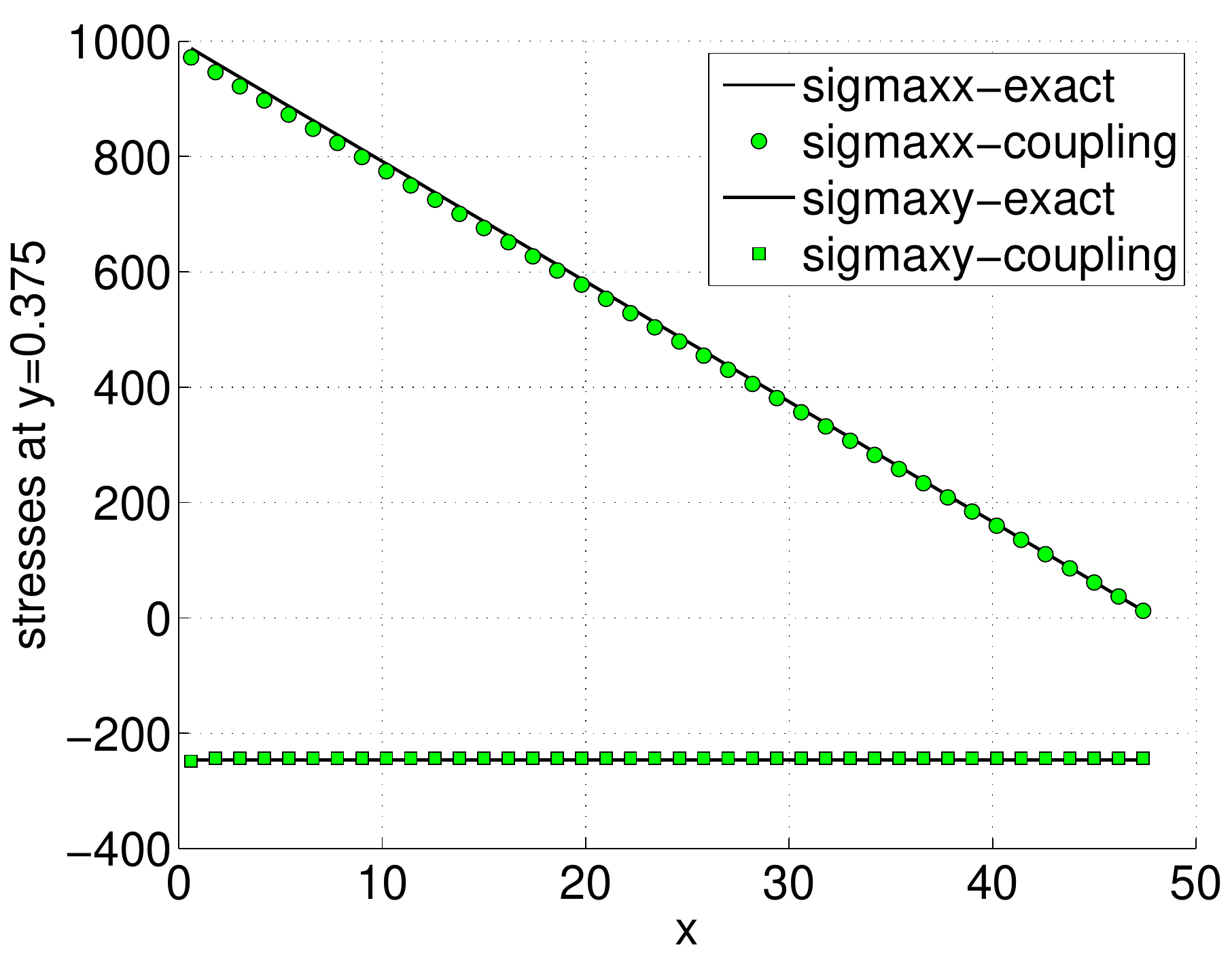}}
   \subfloat[Stresses over the beam height]{\includegraphics[width=0.45\textwidth]{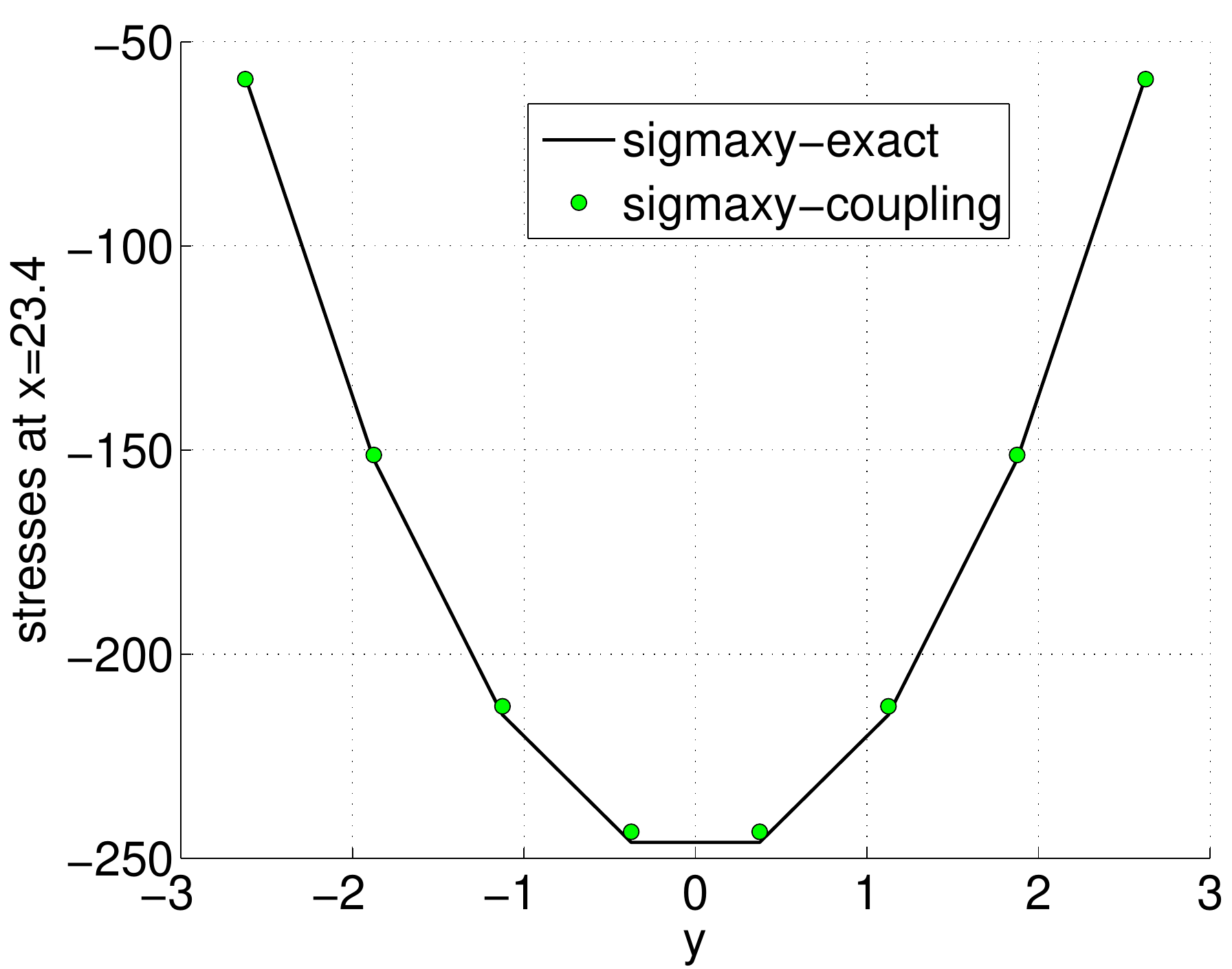}}
  \caption{Timoshenko beam: stresses obtained with a conforming mesh ($20\times8$ for each domain).}
  \label{fig:beam-stress}
\end{figure}

\noindent \textbf{B-splines elements} Next, we study the performance of the B-splines elements of
which one mesh is given in Fig.~\ref{fig:beam-bspline-mesh}. 
Dirichlet BCs are enforced using the least square projection method see \eg \cite{nguyen_iga_review}.
Note that Nitche's method can also be used to weakly enforce the Dirichlet BCs. However, we use
Nitsche's method only to couple the patch interfaces.
As detailed in \cite{hughes-fem-book} for Lagrangian basis functions,
a rule of $(p+1)\times(q+1)$ Gaussian quadrature can be applied for
two-dimensional elements in which $p$ and $q$ denote the orders of
the chosen basis functions in the $\xi$ and $\eta$ direction. The same
procedure is also used for NURBS basis functions in the present work,
although it should be emphasised that Gaussian quadrature is not optimal for IGA
\cite{hughes_efficient_2010,Auricchio201215}. The stresses are given in   Fig.~\ref{fig:beam-bspline-stress}.

\begin{figure}[htbp]
  \centering 
   \includegraphics[width=0.55\textwidth]{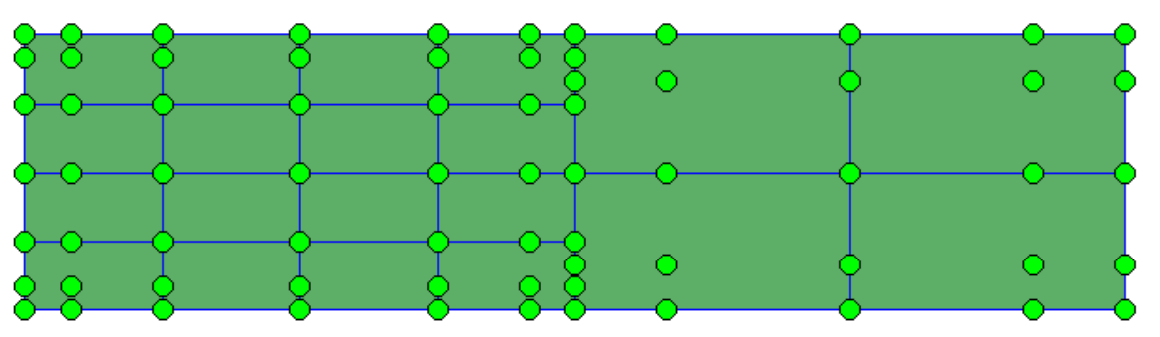}
  \caption{Timoshenko beam: B-spline bi-cubic ($p=q=3$) mesh with $4\times4$ elements for the left
  domain and $2\times2$ elements for the right one. The filled circles denote the control points.}
  \label{fig:beam-bspline-mesh}
\end{figure}

\begin{figure}[htbp]
  \centering 
   \includegraphics[width=0.45\textwidth]{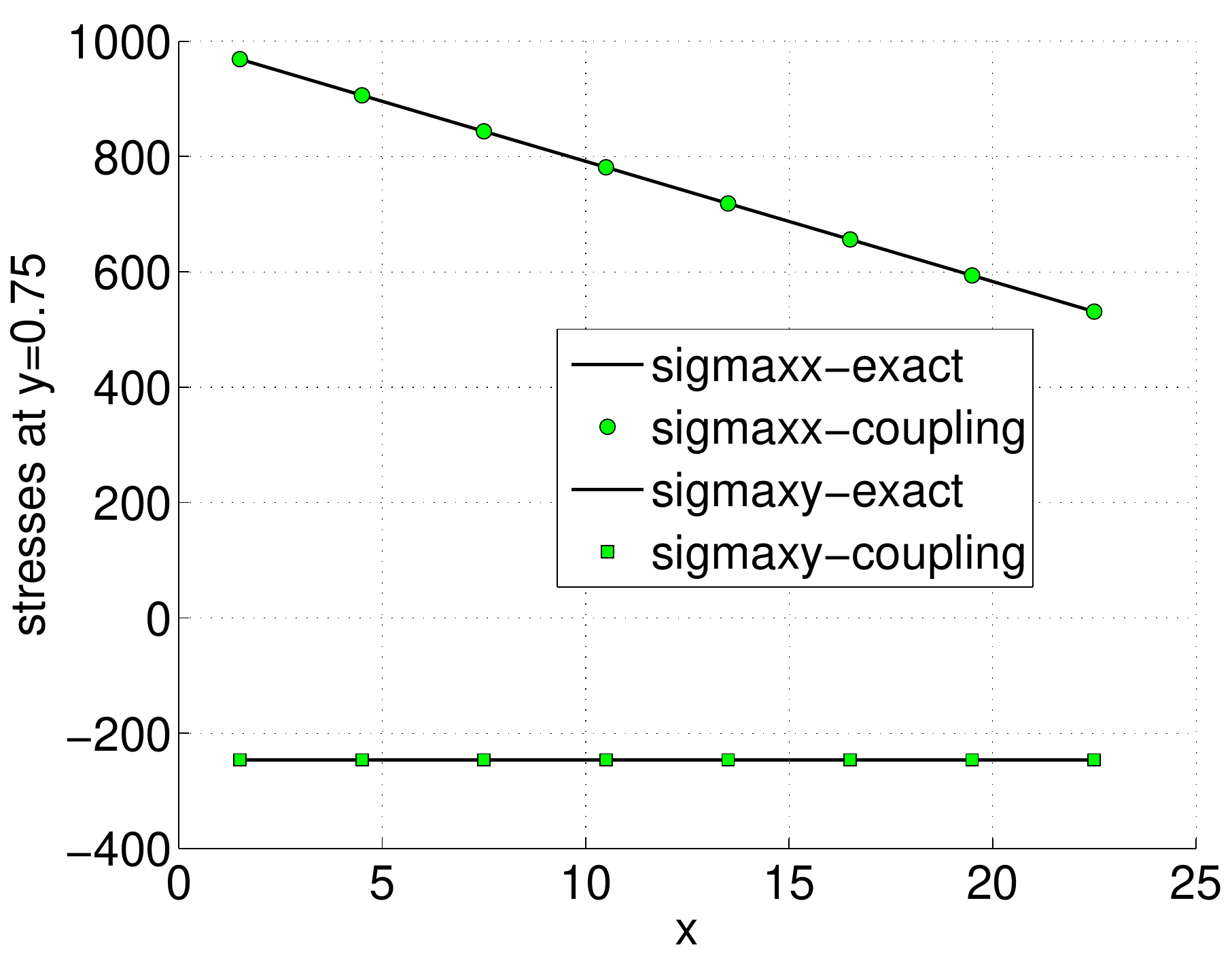}
   \includegraphics[width=0.45\textwidth]{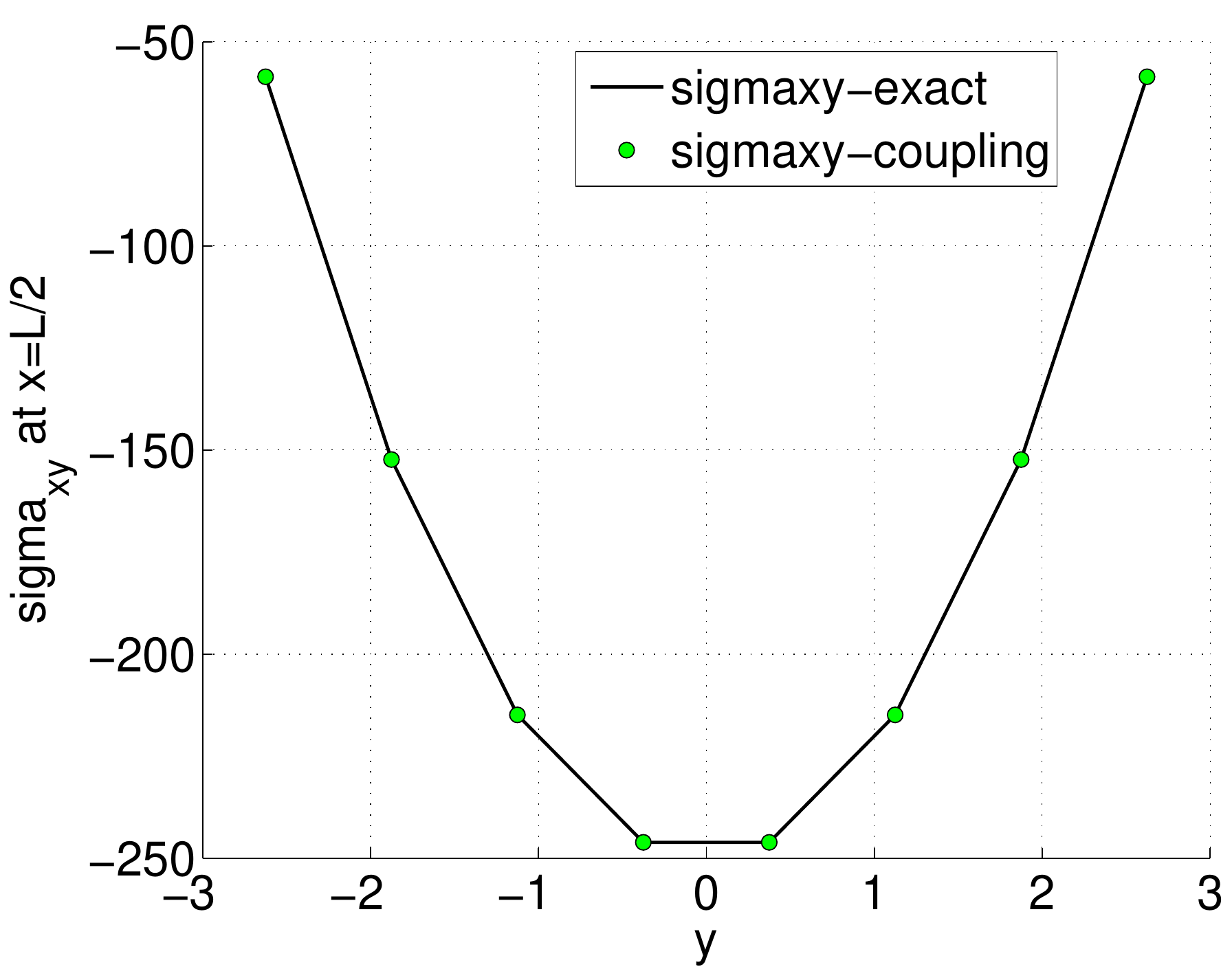}
  \caption{Timoshenko beam: stresses with B-splines elements. The left domain
  is meshed by $8\times8$ cubic elements and the right domain with $2\times2$ cubic elements.}
  \label{fig:beam-bspline-stress}
\end{figure}

Finally we present results obtained with a non-hierarchical B-spline mesh as given in 
Fig.~\ref{fig:beam-bspline-general}: a $8\times6$ bi-cubic mesh is used for the left domain
and a bi-cubic $4\times3$ mesh is used for the right domain. A quadratic stress profile was obtained
where the theoretical maximum value along the midline of the beam (250) can be observed.\\

\begin{figure}[htbp]
  \centering 
   \includegraphics[width=0.55\textwidth]{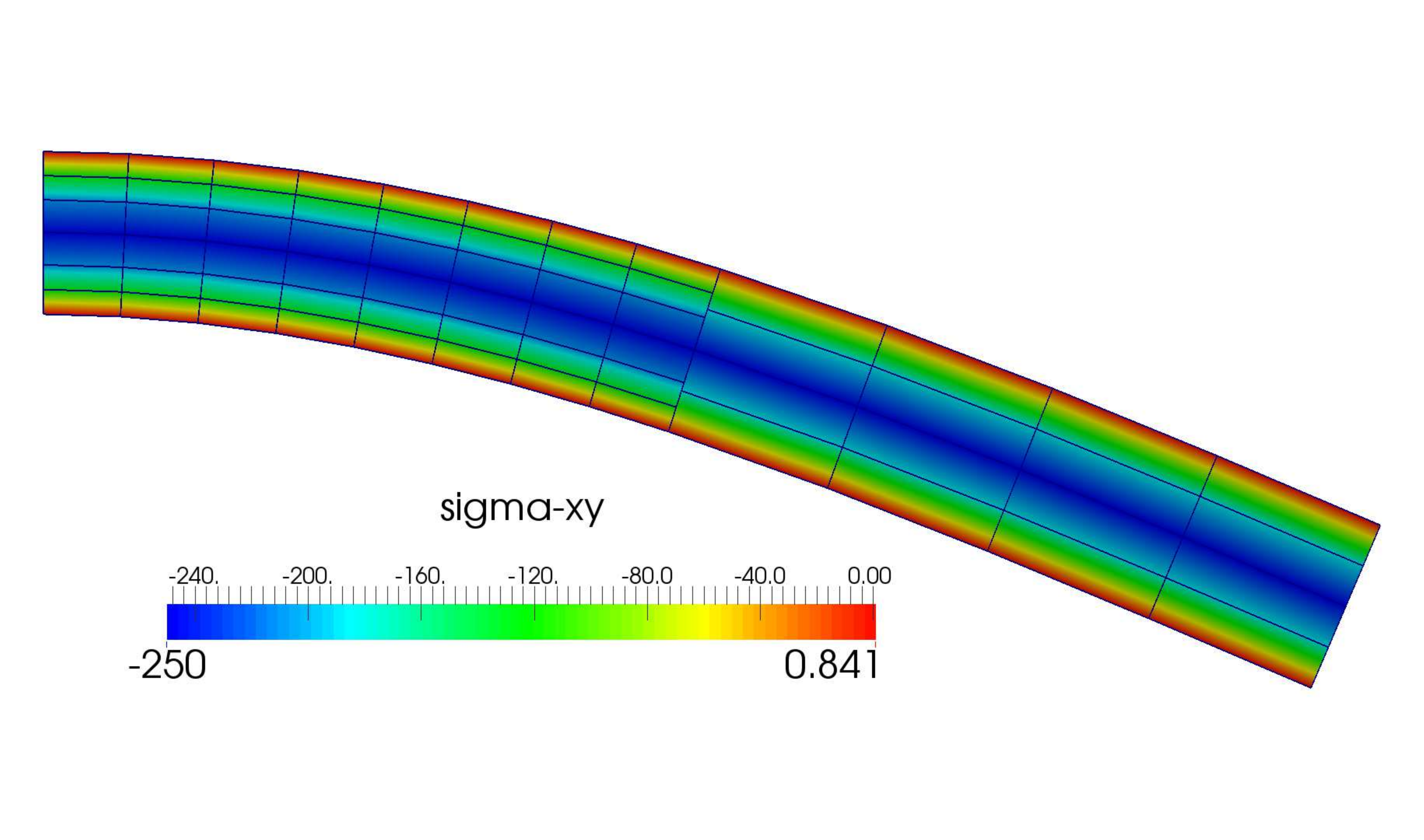}
  \caption{Timoshenko beam: non-hierarchical B-spline mesh ($8\times6$ cubic elements
        for the left domain and $4\times3$ cubic elements for the right domain).}
  \label{fig:beam-bspline-general}
\end{figure}

\noindent \textbf{Convergence study} In order to assess the convergence of the method, 
displacement and energy norms are evaluated with the energy norm given by
\begin{equation}
   e_\text{energy} = \left[\frac{1}{2}\int_{\Omega} \left(\bsym{\varepsilon}_{\mathrm{num}} 
   - \bsym{\varepsilon}_{\mathrm{exact}}\right)\cdot \vm{D} \cdot
   \left(\bsym{\varepsilon}_{\mathrm{num}} - \bsym{\varepsilon}_{\mathrm{exact}} \right) 
   \di \Omega \right]^{\frac{1}{2}},
\end{equation}
and the displacement norm defined as

\begin{equation}
e_\text{displacement} = \left\{{\D \int_{\Omega} \left[ \left( {\bf u}_\text{num} -
{\bf u}_\text{exact} \right) \cdot \left( {\bf u}_\text{num} - {\bf u}_\text{exact}
\right) \right] \di \Omega }\right\}^{1/2},
\end{equation}

\noindent where $\bsym{\varepsilon}_{\mathrm{num}}$, and
$\bsym{\varepsilon}_{\mathrm{exact}}$ are the numerical strain
vector and exact strain vector, respectively. The same
notation applies to the displacement vector ${\bf u}_\text{num}$ and
${\bf u}_\text{exact}$. In the post-processing step, the above norms
are calculated using the same Gauss-Legendre quadrature that has been adopted
for the stiffness matrix computation. 

\begin{figure}[htbp]
  \centering 
   \includegraphics[width=0.55\textwidth]{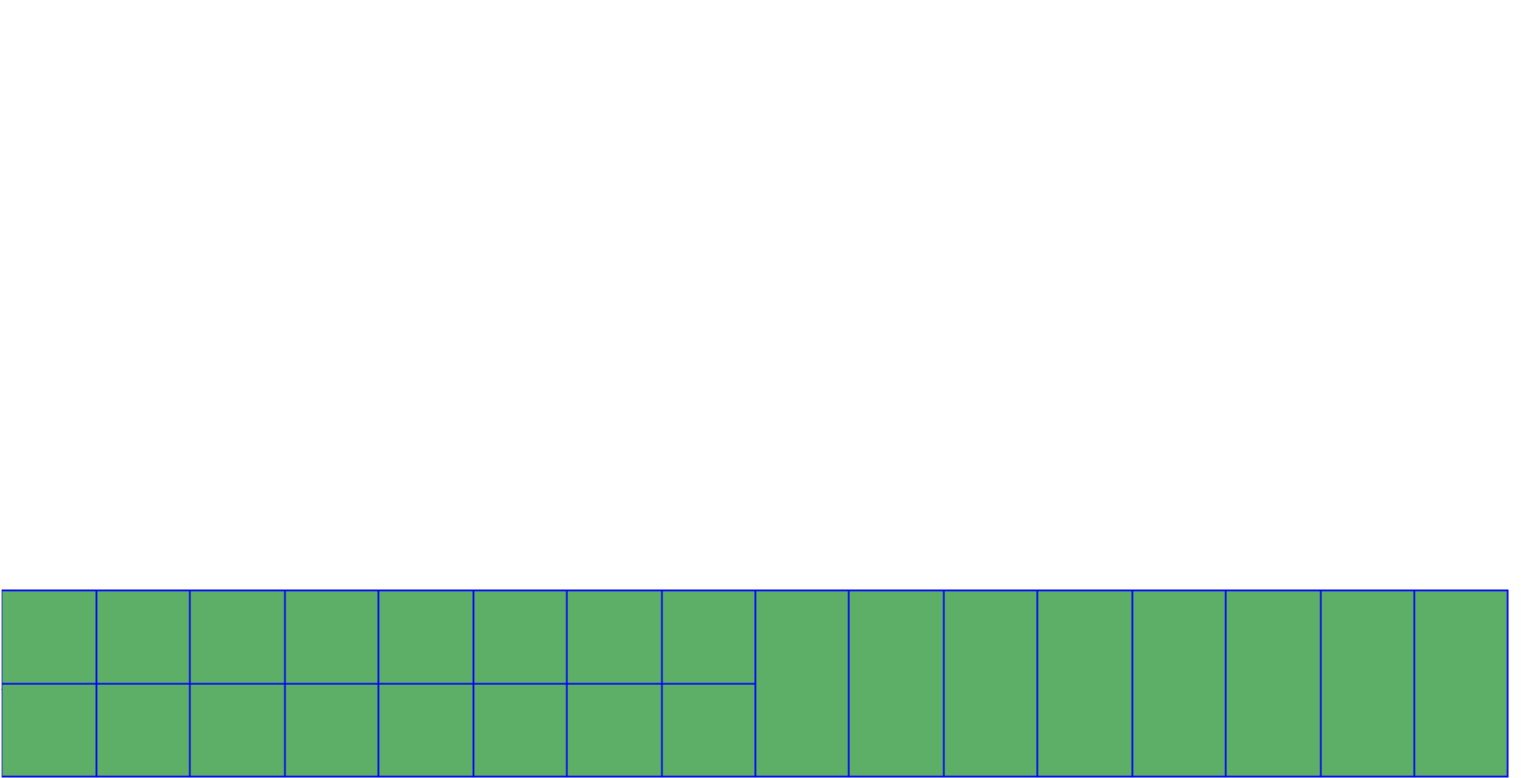}
  \caption{Convergence study of the Timoshenko beam:  initial mesh from which refined meshes are obtained
  by dividing each knot span into two equal halves.}
  \label{fig:beam-convergence-mesh}
\end{figure}

\begin{figure}[htbp]
  \centering 
   \subfloat[displacement norm]{\includegraphics[width=0.49\textwidth]{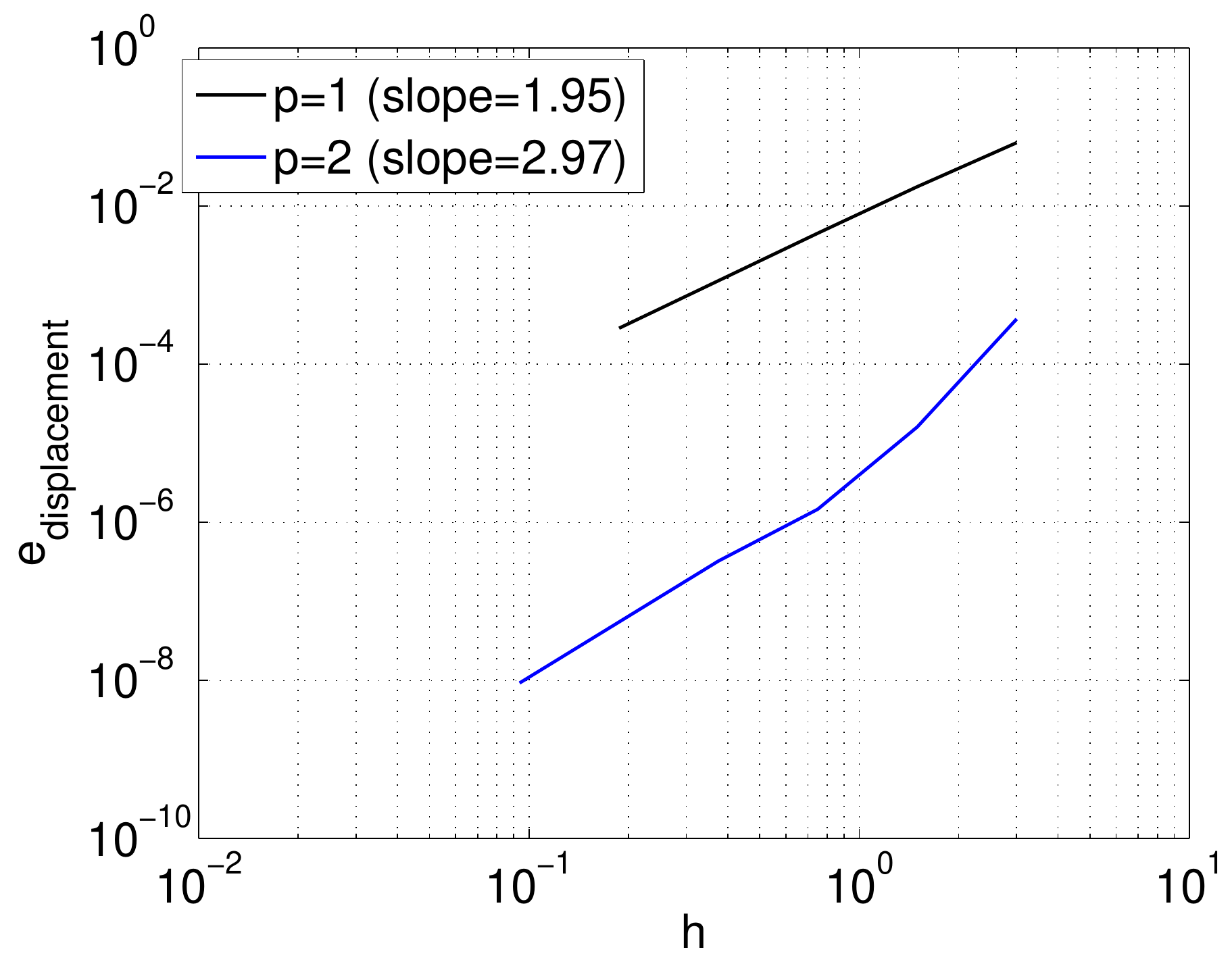}}
   \subfloat[energy norm]{\includegraphics[width=0.49\textwidth]{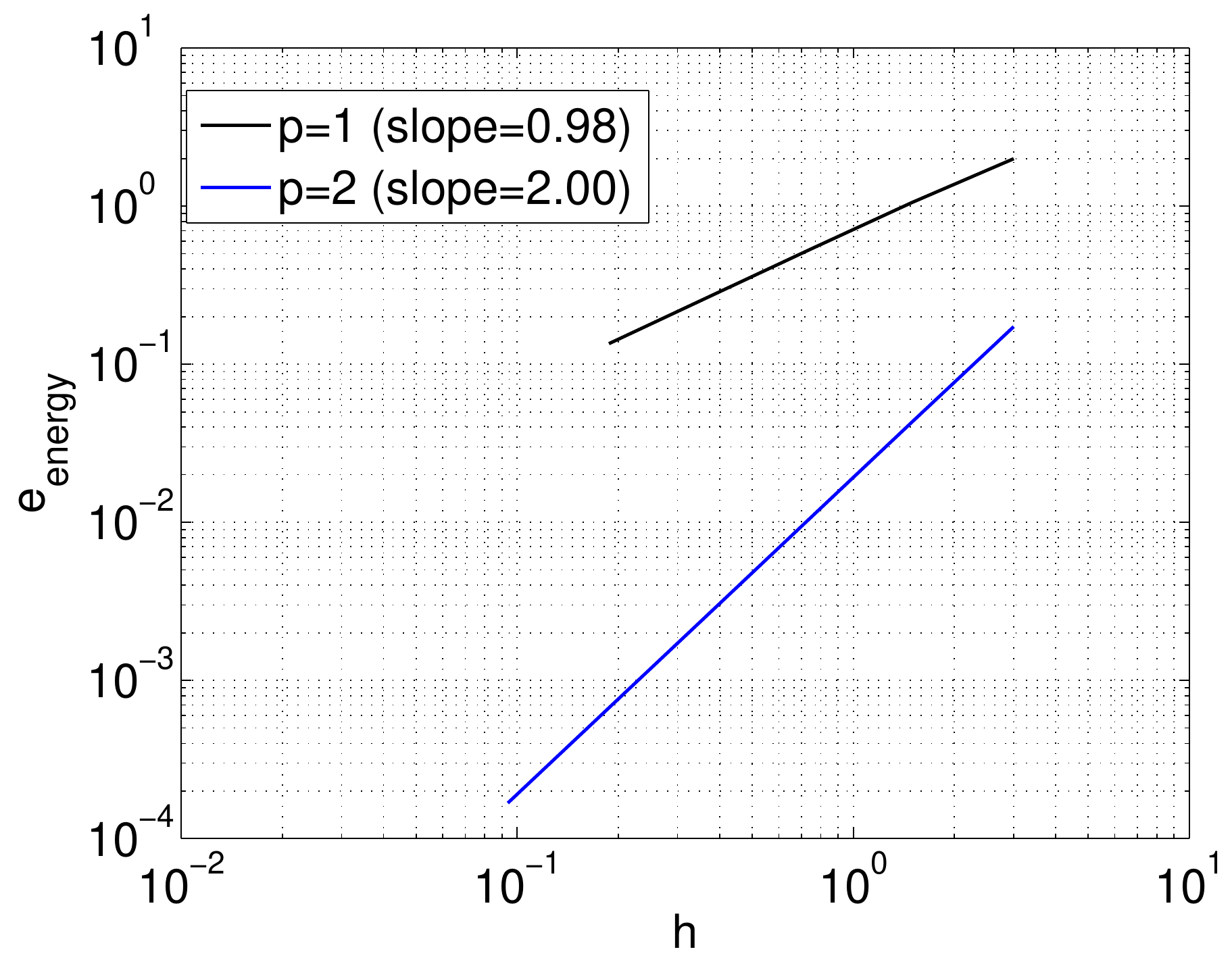}}
  \caption{Timoshenko beam: convergence plots.}
  \label{fig:tbeam-norms}
\end{figure}

The initial mesh from which refined meshes were obtained is given in   Fig.~\ref{fig:beam-convergence-mesh}.
It can be shown that for linear elasticity $\alpha$ depends on the element size $h_e$ and the material parameters, see for example \cite{fritz-nitsche-mortar,Bazilevs200712}

\begin{equation}
\alpha = \frac{\lambda + \mu}{2} \frac{\theta(p)}{h_e}
\label{eq:alpha-estimation}
\end{equation}
where $\theta(p)$ is a positive number that depends only on the polynomial order $p$ of the finite element
approximation. For bilinear basis functions, we set $\theta(p=1)=12$ and for bi-quadratic basis functions, we set
$\theta(p=2)=36$. These values were chosen so that the stiffness matrix is positive definite. Thus, for each mesh,
Equation~\eqref{eq:alpha-estimation} was used to compute the stabilisation parameter.
The convergence plots are given in   Fig.~\ref{fig:tbeam-norms} where optimal convergence rates for both displacement 
and energy norms were obtained.  Note that minimum values of $\alpha$ can be computed based on a 
numerical analysis of the discrete forms and lead to the global \cite{Griebel}
and local generalized eigenvalue approaches \cite{embar_imposing_2010}.

\subsection{Plate with a center inclusion}

Consider a plate with a center inclusion as given in Fig.~\ref{fig:inclusion-geo}. 
The matrix properties are denoted by $E_m$ and $\nu_m$ and the inclusion properties are 
denoted by $E_i$ and $\nu_i$. 
A traction along the vertical direction is applied on the top edge while nodes along the bottom edge 
are constrained. This problem is solved with (1) embedded Nitsche's method and (2) XFEM which are methods
that do not require a mesh conforming to the inclusion. The XFEM mesh is given in Fig.~\ref{fig:inclusion-mesh}a
where $30\times60$ four-noded quadrilateral (Q4) elements are adopted. The material interface is modeled via enrichment
functions (the $abs$ enrichment function) proposed in \cite{Sukumar1}. Meshes in the Nitsche's method, cf. Fig.~\ref{fig:inclusion-mesh}b, consist of
a background mesh for the plate ($32\times64$ Q4 elements) and another mesh for the inclusion which is embedded
in the background mesh ($16\times16$ bi-quadratic NURBS elements).

\begin{figure}[htbp]
  \centering 
   \includegraphics[width=0.45\textwidth]{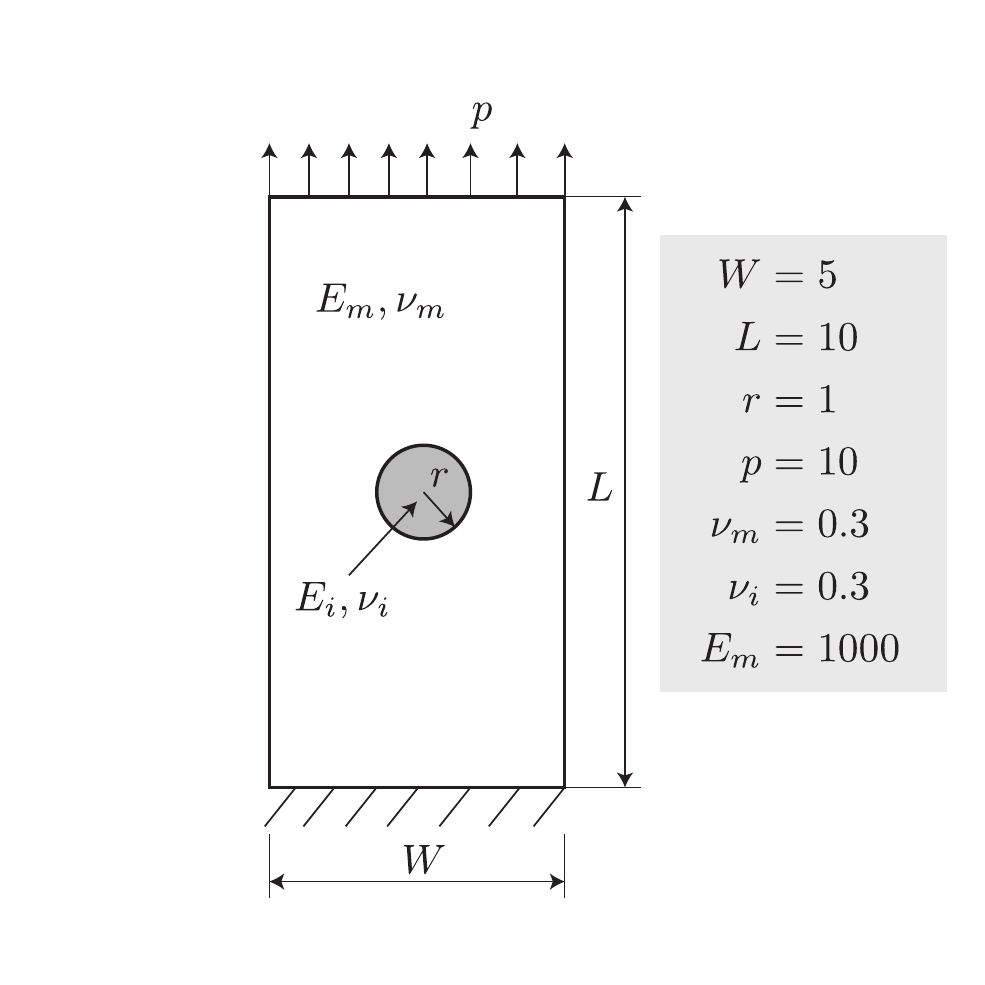}
  \caption{A plate with a center inclusion.}
  \label{fig:inclusion-geo}
\end{figure}

\begin{figure}[htbp]
  \centering 
   \subfloat[XFEM]{\includegraphics[width=0.33\textwidth]{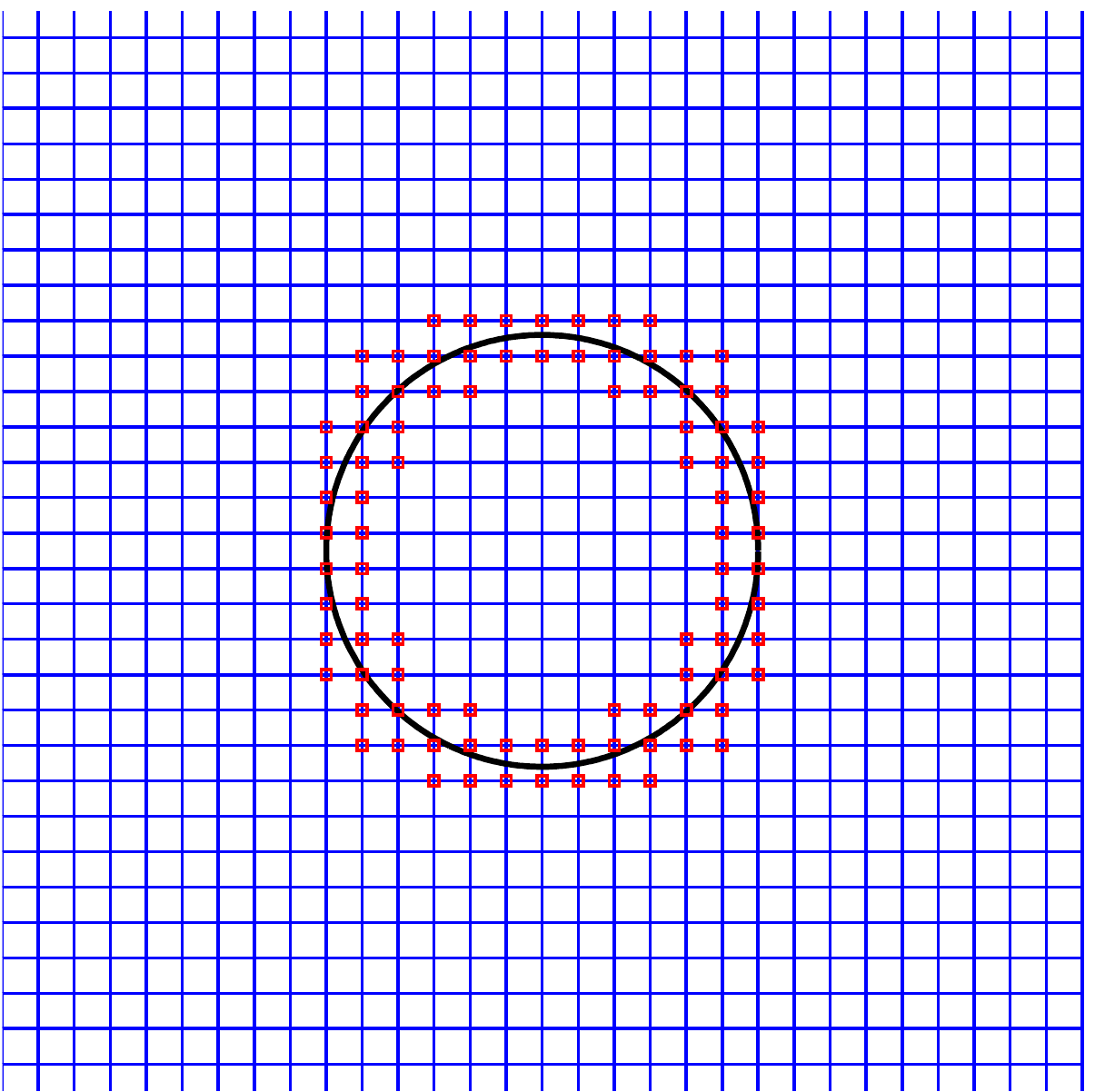}}\;\;
   \subfloat[Nitsche]{\includegraphics[width=0.4\textwidth]{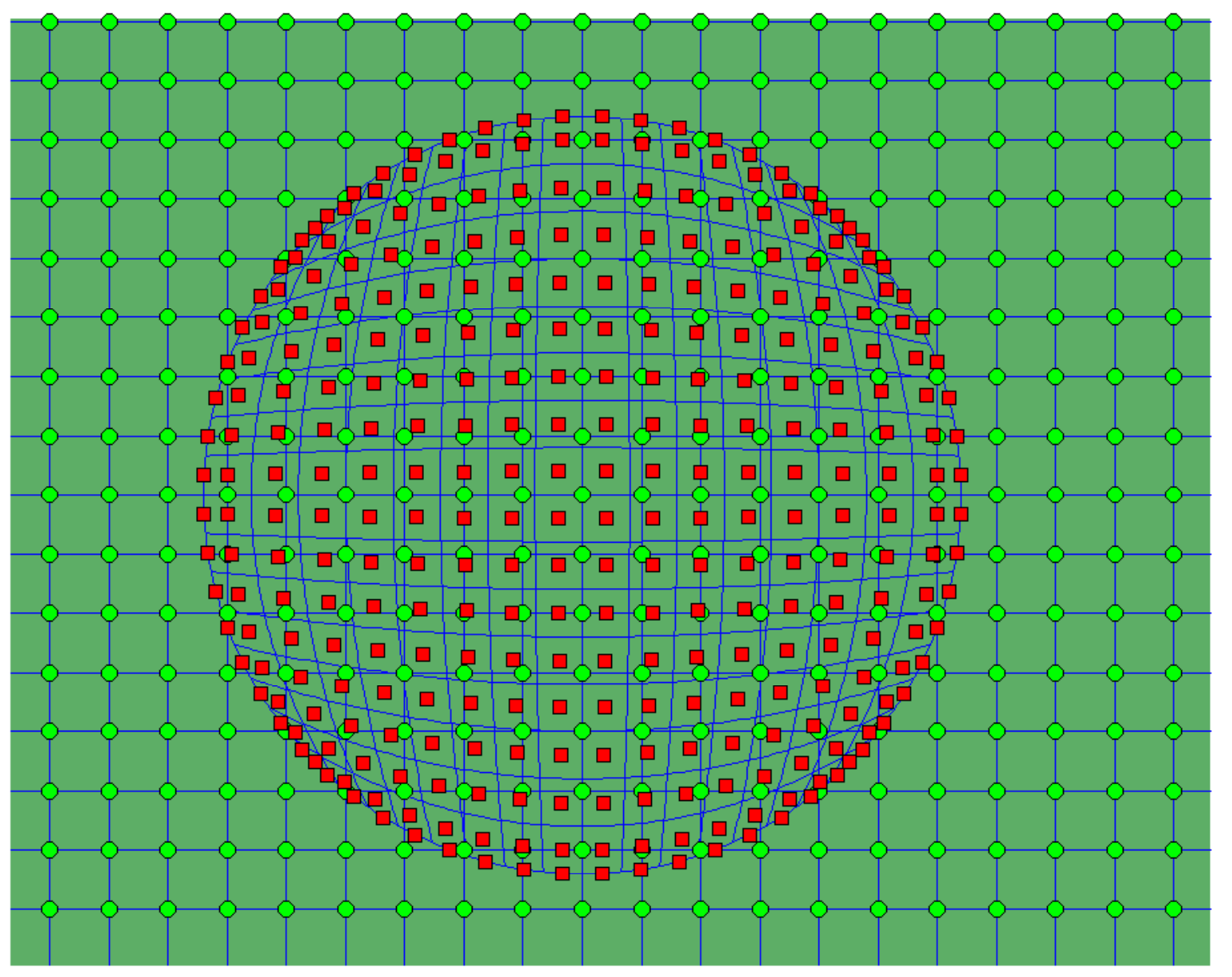}}
  \caption{Plate with a center inclusion: (a) XFEM mesh with enriched nodes and (b) Nitsche's method
  with embedded mesh.}
  \label{fig:inclusion-mesh}
\end{figure}

For details on the Nitsche based embedded mesh method, we refer to \eg \cite{Sanders2011a}. Here, we apply
this method in the context of IGA by using NURBS elements. 
The implementation is briefly explained as follows. The assembly of inclusion elements is standard and
the assembly of background elements is similar to XFEM for voids--void elements (completely covered by
inclusion elements) do not contribute to the total stiffness matrix, cut elements (elements cut by the inclusion)
require special integration scheme in which the part falls within the inclusion domain is not integrated.   
This can be achieved using the standard sub-triangulation technique in the context of XFEM \cite{mos_finite_1999}
or the hierarchical element subdivision employed in the Finite Cell Method \cite{NME:NME4522} or the technique
used in the NEFEM (NURBS Enhanced FEM) \cite{sevillafern'andez-m'endez2008b}. Here, for simplicity, we used
the hierarchical element subdivision method.
We refer to Fig.~\ref{fig:inclusion-nitsche-explain}. 
The inclusion Young's modulus is $E_i=1$.
Due to the contrast in Young's moduli, the average operator given in Equation~\eqref{eq:general-average}
was used with $\gamma=E_m/(E_m+E_i)$ as proposed in \cite{Sanders2011a}. The stabilisation parameter 
is chosen empirically $\alpha=1e6$.
Fig.~\ref{fig:inclusion-nitsche-res} shows the contour plot of $u_y$ solutions obtained with both methods.
A good agreement of Nitsche solution compared with XFEM solution can be observed.

\begin{figure}[htbp]
  \centering 
   \includegraphics[width=0.45\textwidth]{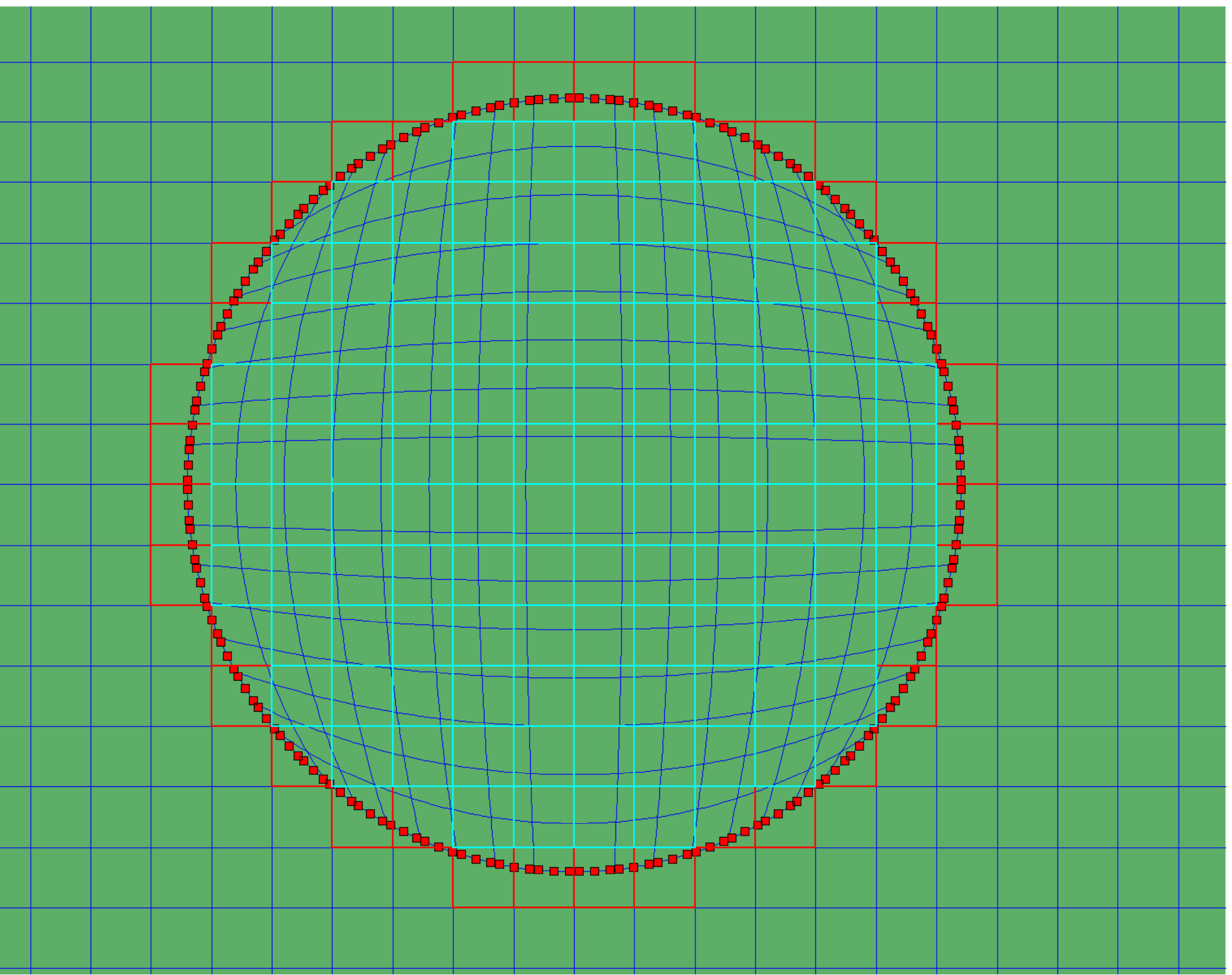}
  \caption{A plate with a center inclusion: Nitsche based embedded mesh method.
  The red filled squares denote Gauss points to evaluate the coupling matrices. Cyan squares denote void elements
  and red squares represent cut elements.}
  \label{fig:inclusion-nitsche-explain}
\end{figure}

\begin{figure}[htbp]
  \centering 
   \includegraphics[width=0.45\textwidth]{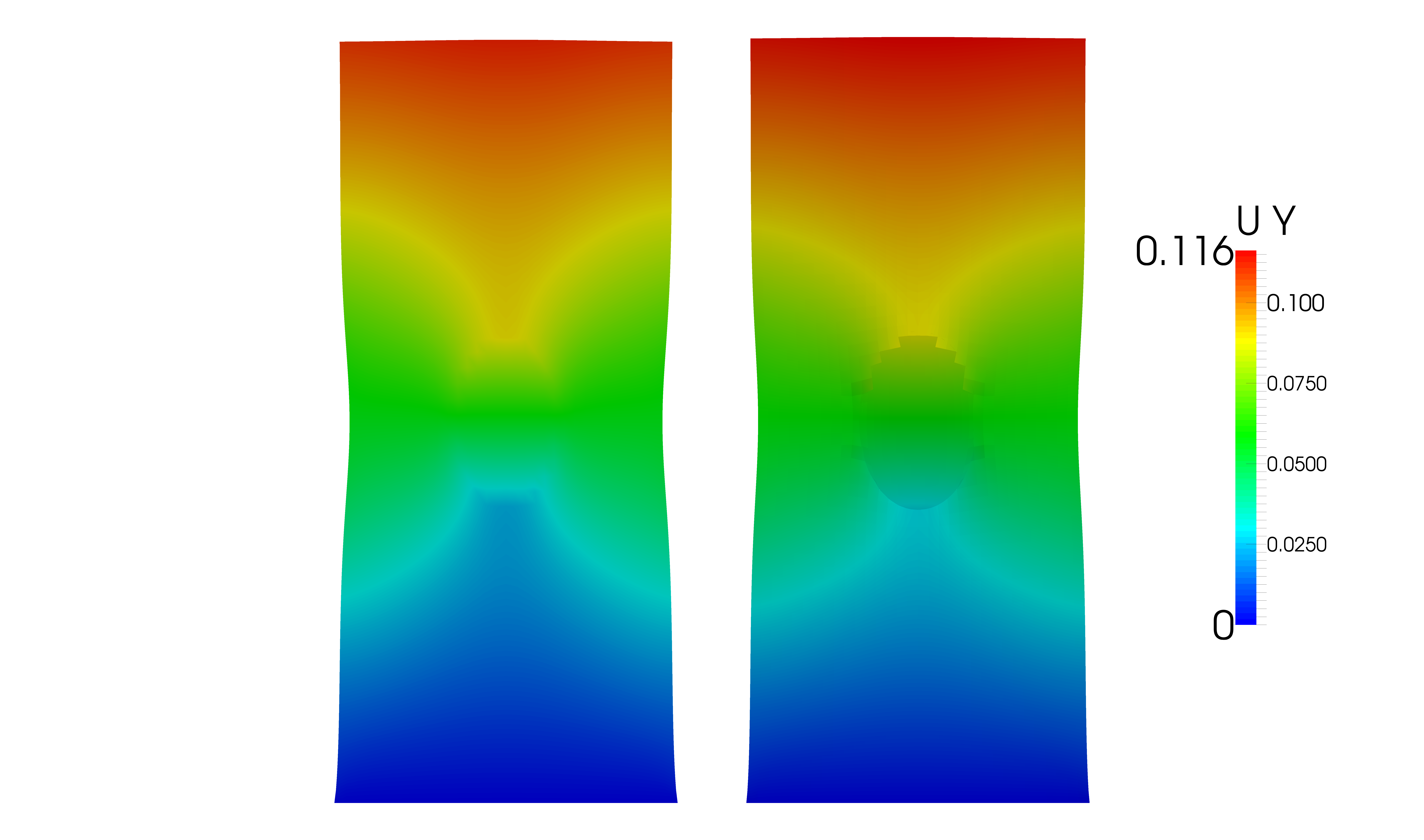}
  \caption{A plate with a center inclusion: contour plot of $u_y$ solutions--xfem (left) and Nitsche (right).}
  \label{fig:inclusion-nitsche-res}
\end{figure}

\subsection{3D-3D coupling}

In order to test the implementation for 3D problems, we consider the 3D cantilever beam shown in 
Fig.~\ref{fig:beam3D-geo}.
The data are: Young's modulus $E=1000$, Poisson's ratio $\nu=0.3$, $L=10$, $W=H=1$ and the imposed
displacement in the $z$-direction is $1$. The non-conforming 
B-splines discretisation is given in Fig.~\ref{fig:beam3D-mesh} where the beam is divided into two
equal parts. A value of ... was used for the stabilisation parameter $\alpha$.
In Fig.~\ref{fig:beam3D-res} the contour plot of $\sigma_{xx}$ is given and a comparison was made with
a standard Galerkin discretisation of $32\times4\times4$ tri-cubic B-splines elements and a good agreement was obtained.

\begin{figure}[htbp]
  \centering 
   \includegraphics[width=0.65\textwidth]{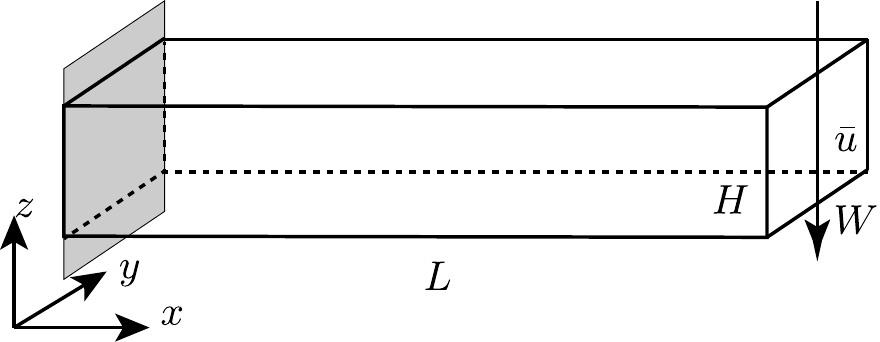}
  \caption{A 3D cantilever beam subjected to an imposed vertical displacement.}
  \label{fig:beam3D-geo}
\end{figure}

\begin{figure}[htbp]
  \centering 
   \includegraphics[width=0.65\textwidth]{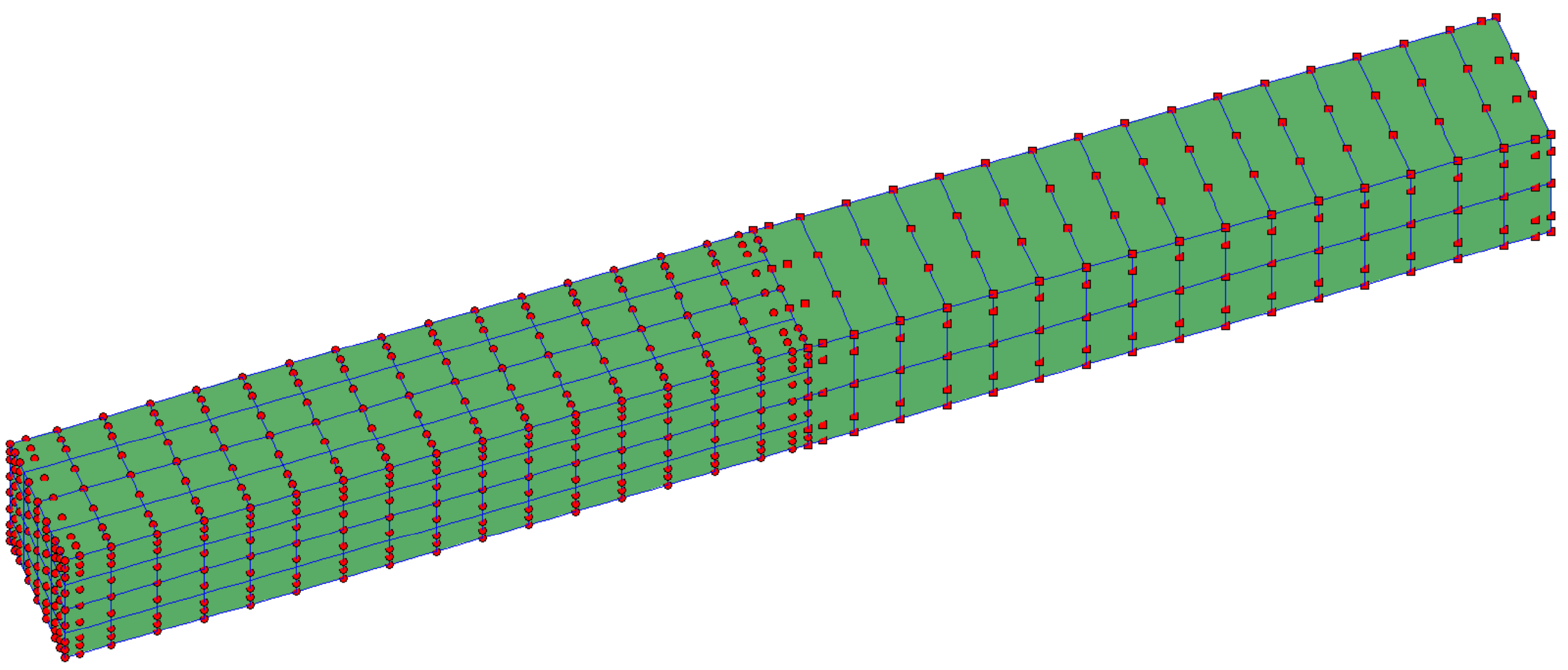}
  \caption{A 3D cantilever beam subjected to an imposed vertical displacement: $16\times4\times4$ tri-cubic
     B-splines elements for the left domain and  $16\times1\times2$ tri-cubic
     B-splines elements for the right domain.}
  \label{fig:beam3D-mesh}
\end{figure}

\begin{figure}[htbp]
  \centering 
   \subfloat[Nitsche]{\includegraphics[width=0.55\textwidth]{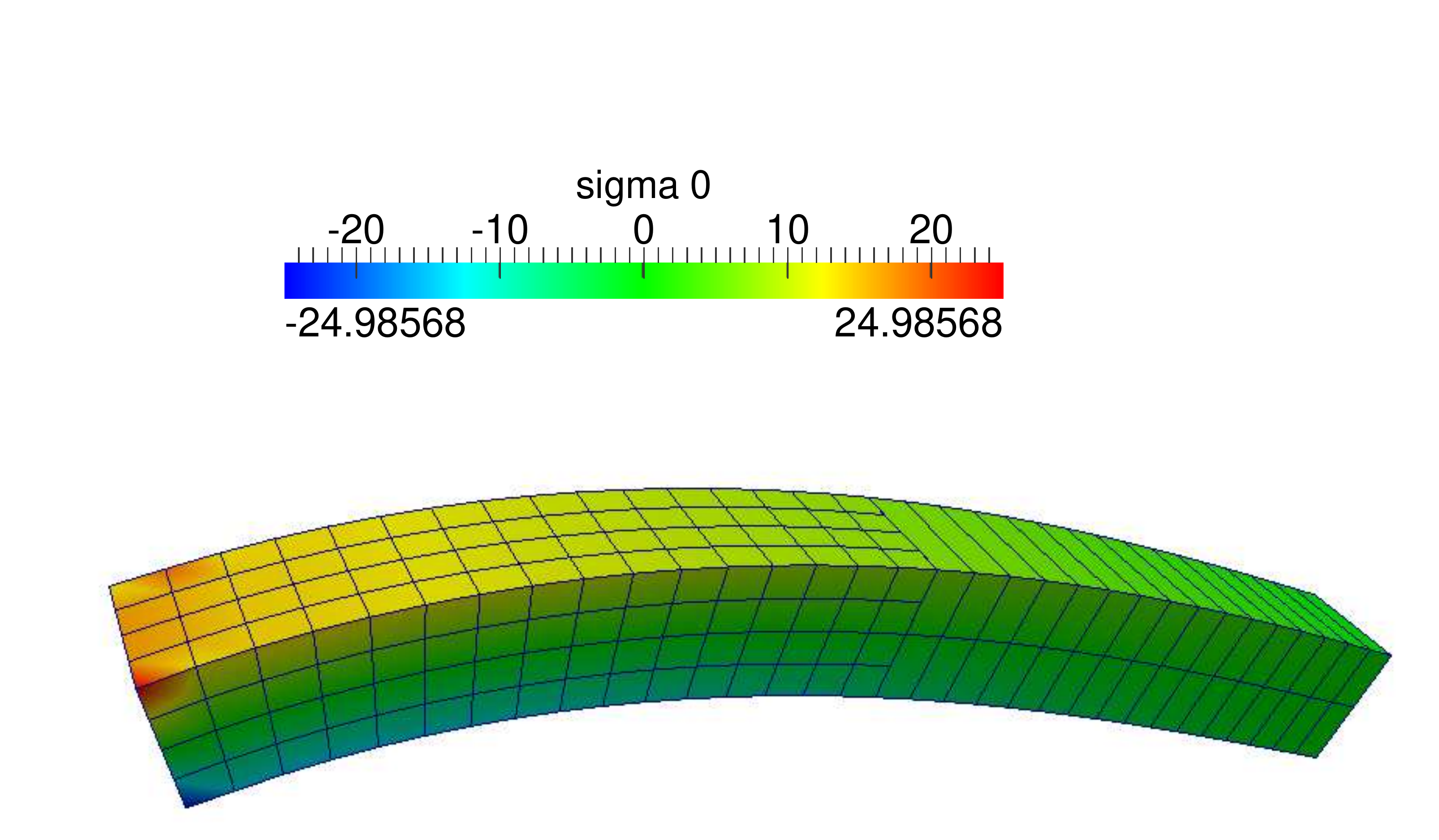}}\\
   \subfloat[Galerkin]{\includegraphics[width=0.55\textwidth]{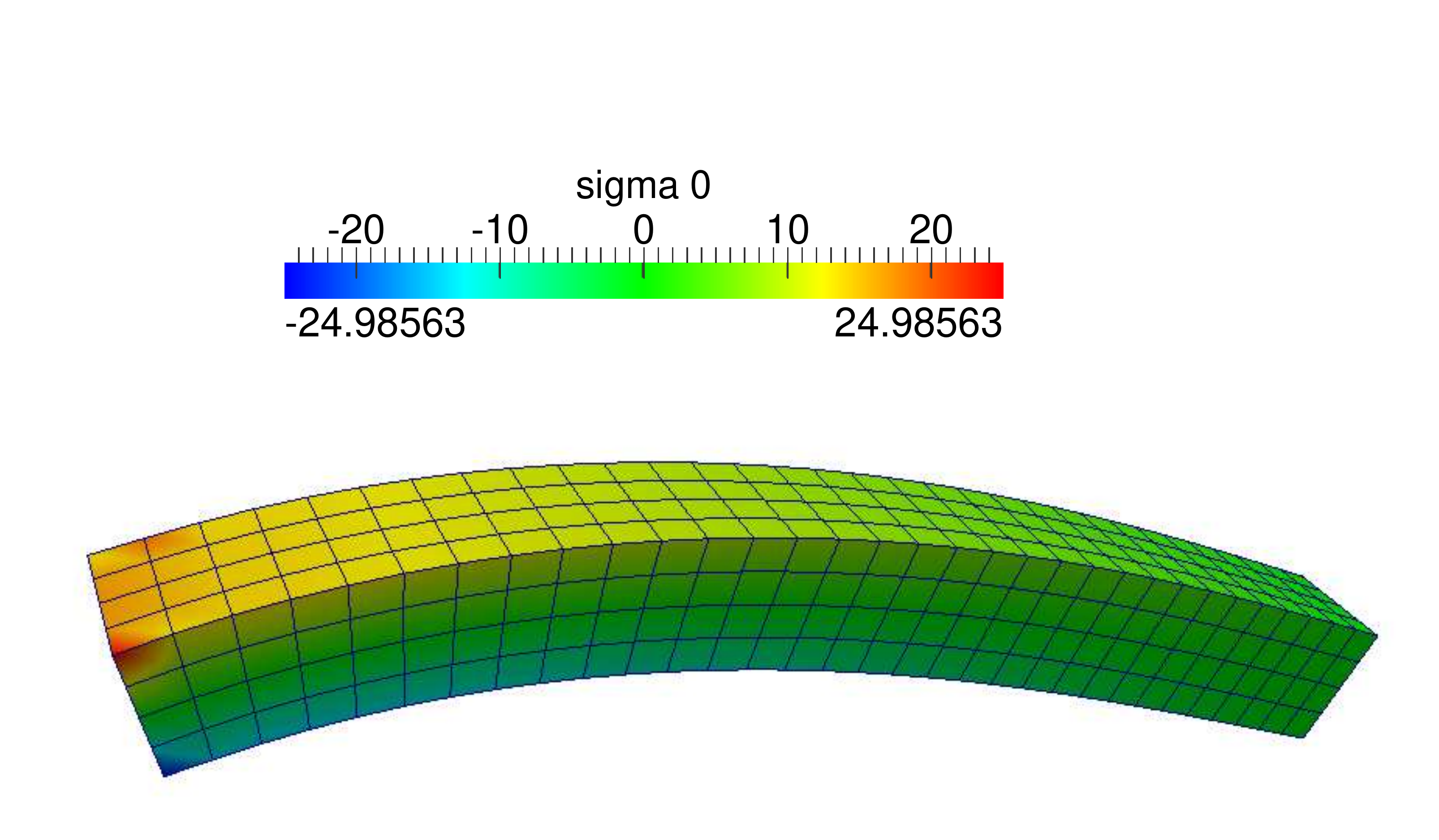}}
  \caption{Timoshenko beam.}
  \label{fig:beam3D-res}
\end{figure}

\subsection{Connecting rod}

The method is now applied to a more complicated geometry, taking into account more than one interface coupling, curved interfaces and interfaces with different dimension.
This geometry is a simplified representation of a \emph{connecting rod}, which is a component of an internal combustion engine, and represents a classic linear case in the stress-strain static analysis.
The geometric input model is composed by three NURBS patches (see Fig.~\ref{fig:concepts}) with two coupling interfaces. 
The dimensions are consistent with an actual component and the material properties are Young's modulus $E=2\times10^5$ MPa, Poisson's ratio $\nu=0.3$  which come from a standard steel material.
Boundary conditions are represented in Fig.~\ref{fig:conrod3D-BC}: ideal fixed boundary condition on the two vertical surfaces of the (\emph{big-end}) and a vertical total force $F=1000$ N load applied to the internal ring of the \emph{small-end}, according to the effect of the \textit{pin-piston} sub-assembly that transmits a bending moment to the connecting-rod stem.
For the simulation the model is refined with tri-cubic functions and $32\times4\times8$ elements for patch 1, $24\times12\times4$ elements for patch 2 and $64\times4\times8$ elements for patch 3, resulting in a total number of 4224 elements and 11305 control points.
For both coupling interfaces the smaller faces are the regions where the surface integration is performed and a stabilization parameter $\alpha=1\times10^8$ was chosen empirically.
The results are shown in Fig.~\ref{fig:conrod3D-results}, where displacement and stress fields are plotted.
The displacement distribution is the typical progressive cubic polynomial form of the analytical Saint-Venant model.
The pattern distribution of the Von Mises equivalent failure criterion is used for the comparison of the simulation results in IGA approach with respect to \textit{Siemens-NX} (traditional FE model, discretized with second order tetrahedra, 6182 elements and 11332 nodes Fig.~\ref{fig:NASTRAN-results}).
Typical combined compressive and bending stress/action of the connecting-rod stem is representable with Von Mises stresses closed to zero in the mean plane; superior fibres has the maximum value of traction symmetrically equivalent to the compression of inferior fibres, due to the strictly positive equivalent measure of Von Mises yield criterion.
In both analyses interesting three-dimensional effects are detected: maximum stress values correspond to the free fibres of the stem in superior and inferior surfaces that interact with the big-end; the interaction between the stem and both the big-end and small-end produces an increasing stress value in the azure region in proximity of the neutral axis that is very well described in both analysis, thus demonstrating the IGA model effectiveness of the links between patches; the boundary conditions are typically hyperstatic and only the inner part of the big-end transmits traction/compression reactions (green regions); due to this particular load case, parts of the big-end (blue regions) are superfluous in both analyses and could be deleted, reducing the mass of the component; the internal stress distribution in the inner ring of the small-end shows again very good agreement of the combined compressive and bending stress/action behaviour that reaches the pin region.

\begin{figure}[htbp]
  \centering 
   \includegraphics[width=.7\textwidth]{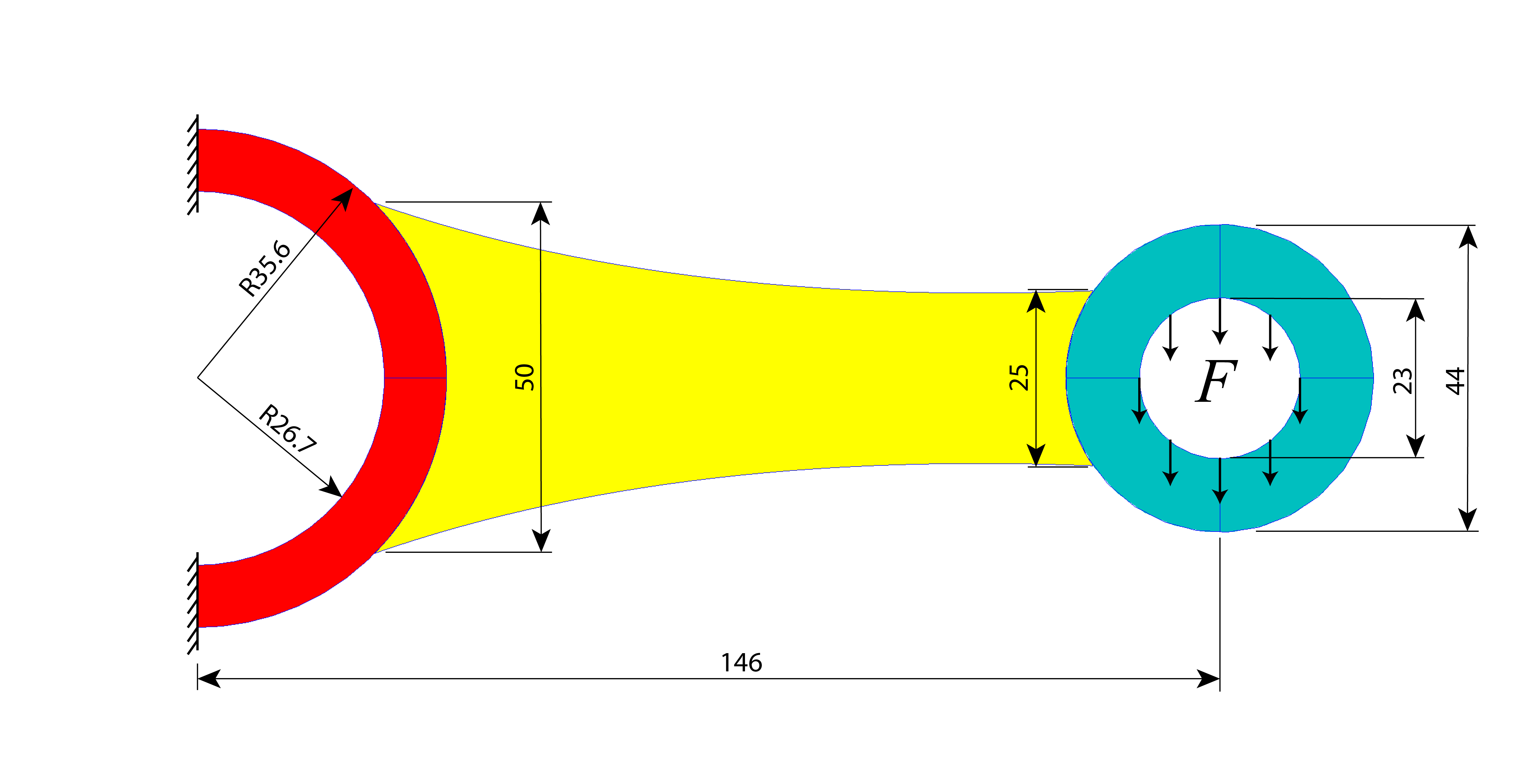}
  \caption{Connecting-rod: geometry and boundary conditions. The dimensions are in mm.}
  \label{fig:conrod3D-BC}
\end{figure}

\begin{figure}[htbp]
  \centering 
   \subfloat[z-displacement field]{\includegraphics[width=.7\textwidth]{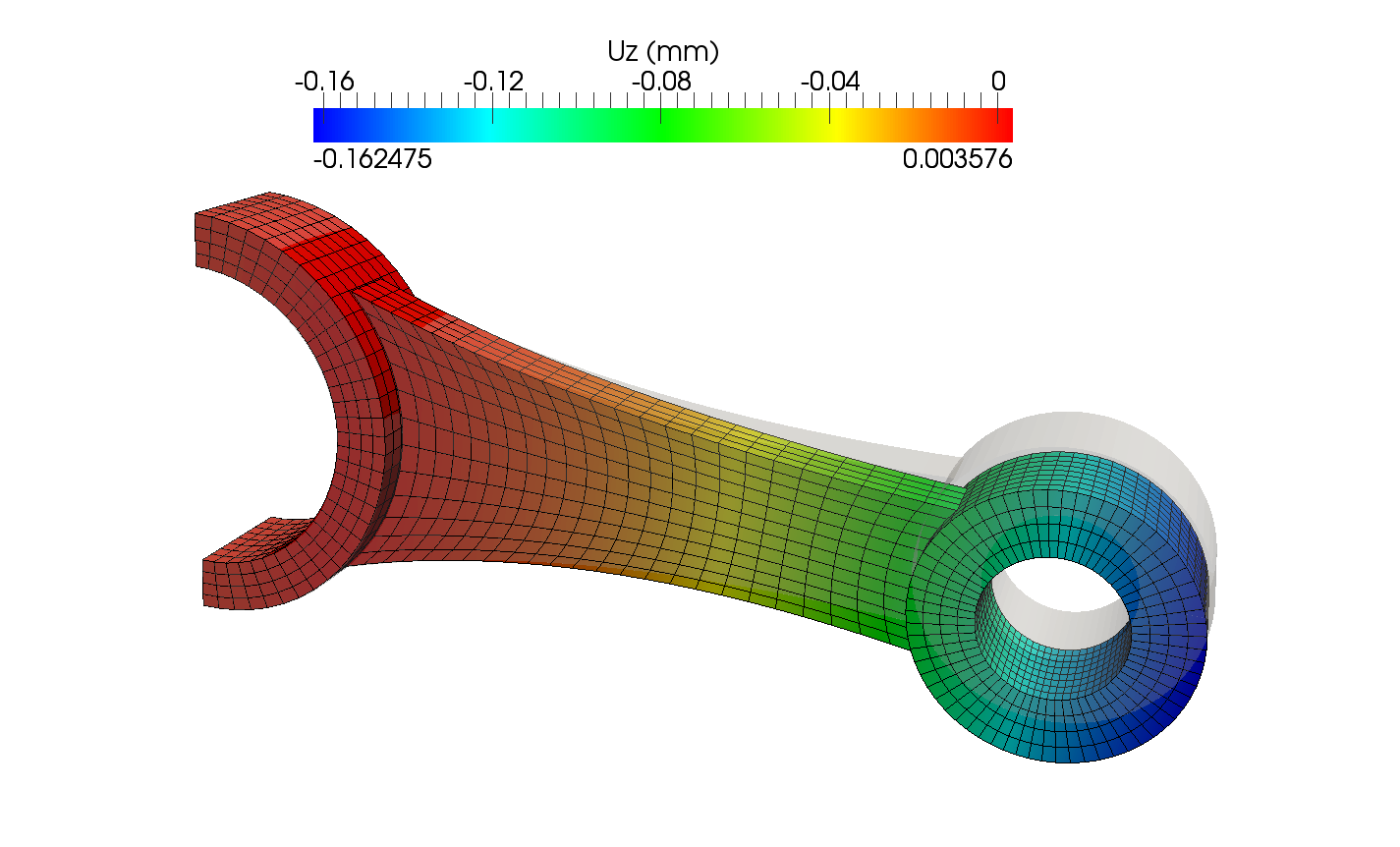}}\\
   \subfloat[Stress field]{\includegraphics[width=.7\textwidth]{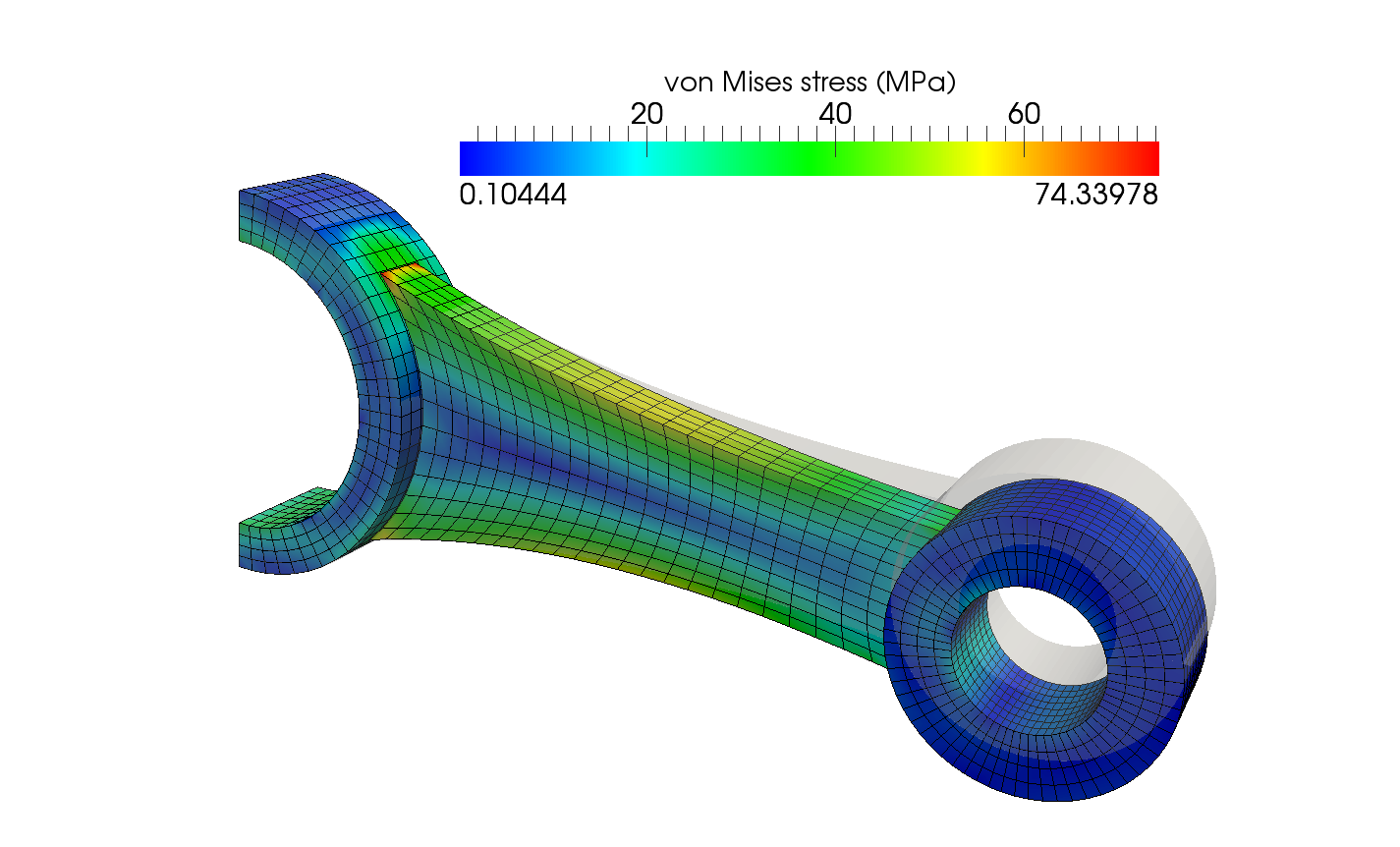}}
  \caption{Results of the connecting rod.}
  \label{fig:conrod3D-results}
\end{figure}

\begin{figure}[htbp]
  \centering 
   \includegraphics[width=1\textwidth]{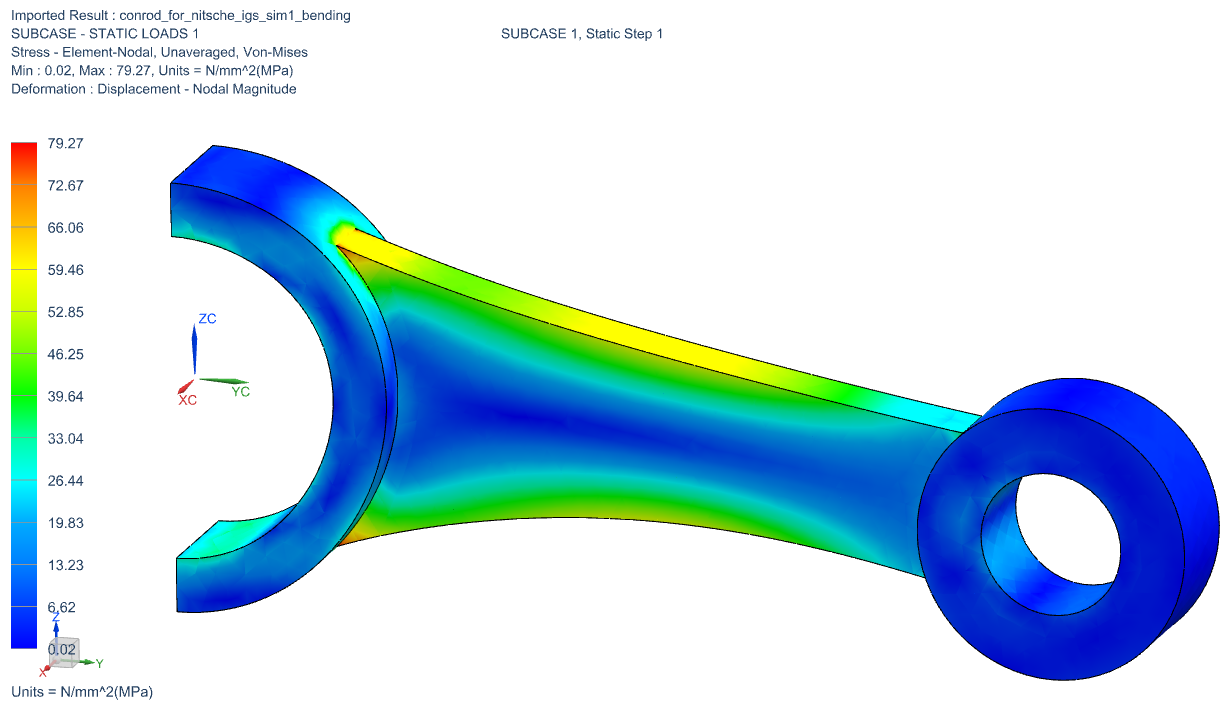}
  \caption{Stress plot from the commercial code \textit{NX-NASTRAN}.}
  \label{fig:NASTRAN-results}
\end{figure}

\section{Conclusions}\label{sec:conclusions}

We presented a Nitsche's method to couple non-conforming NURBS patches. 
Detailed implementation was provided and numerical examples demonstrated the good performance of the method.
The proposed method certainly enlarges the applicability of NURBS based isogeometric analysis.

The contribution was limited to linear elastostatic problems and extension of the method to (1) dynamics problems
and (2) nonlinear material problems is under investigation before one could claim whether Nitsche coupling would be
a viable method for multi-patch NURBS based isogeometric analysis.

As we were preparing the paper for submission, we became aware of contemporary work had been presented the previous week at the US National Congress for Computational Mechanics \cite{Dominik_Nitsche3D} in the context of the finite cell method.

\section*{Acknowledgements}

The authors would like to acknowledge the partial financial support of the
Framework Programme 7 Initial Training Network Funding under grant number
289361 ``Integrating Numerical Simulation and Geometric Design Technology".
St\'{e}phane Bordas also thanks partial funding for his time provided by
1) the EPSRC under grant EP/G042705/1 Increased Reliability for Industrially
Relevant Automatic Crack Growth Simulation with the eXtended Finite Element
Method and 2) the European Research Council Starting Independent Research
Grant (ERC Stg grant agreement No. 279578) entitled ``Towards real time multiscale simulation of cutting in
non-linear materials with applications to surgical simulation and computer
guided surgery''. Marco Brino thanks Politecnico di Torino for the funding that supports his visitor to
iMAM at Cardiff University.

\bibliography{isogeometric}
\bibliographystyle{unsrt}


\end{document}